\documentclass{amsart}
\usepackage{math,amsrefs}
\usepackage[all]{xy}
\newif\ifbilbo
\bilbofalse

\newcommand\Hdim{{\operatorname{Hdim}}}
\newcommand\GKdim{{\operatorname{GKdim}}}
\newcommand\mat[4]{\begin{pmatrix}#1 & #2\\ #3 & #4\end{pmatrix}}
\newcommand\smallmat[4]{(\begin{smallmatrix}#1 & #2\\ #3 & #4\end{smallmatrix})}
\newcommand\vmat[2]{(\begin{smallmatrix}#1\\ #2\end{smallmatrix})}
\newcommand\hmat[2]{(\begin{smallmatrix}#1 & #2\end{smallmatrix})}
\renewcommand\AA{{\mathfrak A}}
\newcommand\DD{{\mathfrak D}}
\newcommand\EE{{\mathfrak E}}
\newcommand\FF{{\mathfrak F}}
\newcommand\II{{\mathfrak I}}
\newcommand\JJ{{\mathfrak J}}
\newcommand\KK{{\mathfrak K}}
\newcommand\PP{{\mathfrak P}}
\newcommand\gr{{\operatorname{\mathfrak{gr}\,}}}
\newcommand\unil{{\operatorname{nil*}}}
\newcommand\jac{{\operatorname{rad}}}

\renewcommand\Gg{{\mathsf{G}}}

\begin{document}
\title{Branch Rings, Thinned Rings, Tree Enveloping Rings}
\author{Laurent Bartholdi}
\address{Laurent Bartholdi, \'Ecole Polytechnique F\'ed\'erale de Lausanne, Institut de G\'eom\'etrie, Alg\`ebre et Topologie, 1015 Lausanne, Switzerland}
\email{laurent.bartholdi@epfl.ch}
\thanks{The author acknowledges support from TU Graz and UC Berkeley, where part of this research was conducted.}
\date{First version November 11, 2004; corrected January 24, 2005}
\begin{abstract}
  We develop the theory of ``branch algebras'', which are
  infinite-dimensional associative algebras that are isomorphic, up to
  taking subrings of finite codimension, to a matrix ring over
  themselves. The main examples come from groups acting on trees.

  In particular, for every field $\Bbbk$ we construct a
  $\Bbbk$-algebra $\KK$ which
  \begin{itemize}
  \item is finitely generated and infinite-dimensional, but has only
    finite-dimensional quotients;
  \item has a subalgebra of finite codimension, isomorphic to
    $M_2(\KK)$;
  \item is prime;
  \item has quadratic growth, and therefore Gelfand-Kirillov dimension
    $2$;
  \item is recursively presented;
  \item satisfies no identity;
  \item contains a transcendental, invertible element;
  \item is semiprimitive if $\Bbbk$ has characteristic $\neq2$;
  \item is graded if $\Bbbk$ has characteristic $2$;
  \item is primitive if $\Bbbk$ is a non-algebraic extension of $\F[2]$;
  \item is graded nil and Jacobson radical if $\Bbbk$ is an algebraic
    extension of $\F[2]$.
  \end{itemize}
\end{abstract}
\maketitle

\section{Introduction}

Although rings arising from groups are very interesting from a ring
theorists' perspective, they are in a sense ``too large'', because
some proper quotient of them may still contain a copy of the original
group. The process of ``quotienting out extra material'' from a group
ring while retaining the original group intact is the ``thinning
process'' described in~\cite{sidki:primitive}.

In this paper, we consider a natural ring arising from a group acting
on a rooted tree, which we call its ``tree enveloping ring''. This is
a re-expression, in terms of matrices, of Said Sidki's
construction~\cite{sidki:primitive}. If the group's action has some
self-similarity modeled on the tree's self-similarity, we may expect
the same to happen for the associated ring, and we use this
self-similarity as a \emph{leitmotiv} for all our results.

Loosely speaking (See \S\ref{ss:alg:branch} for a more precise
statement), a weakly branch algebra is an algebra $\AA$ such that (1)
there is an embedding $\psi:\AA\to M_d(\AA)$ for some $d$, and (2) for
any $n$ there is an element of $\AA$ such that $\psi^n(a)$ has a
single non-zero entry. We show (Theorem~\ref{thm:nopi}) that such
algebras may not satisfy a polynomial identity.

The main construction of weakly branch algebras is via groups acting on trees;
the algebra $\AA$ is then the linear envelope of the groups' linear
representation on the boundary of the tree. We show (Theorem~\ref{thm:gkdim})
that if the groups' orbits on the boundary have polynomial growth of degree
$d$, then the Gelfand-Kirillov dimension of $\AA$ is at most $2d$. In
particular contracting groups generate algebras of finite Gelfand-Kirillov
dimension.

We next concentrate in more detail on the rings $\AA$ arising from the
group $\Gg$ introduced by Grigorchuk in~\cite{grigorchuk:burnside}.
Recall that $\Gg$ is a just-infinite, finitely generated torsion
group.  The algebra $\AA$ over the field $\F[2]$ was already studied
by Ana Cristina Vieira in~\cite{vieira:modular}. The following theorem
summarizes our results in relation with $\Gg$:
\begin{thm}
  The ring $\AA$ is just-infinite and prime
  (Theorem~\ref{thm:galg:jip}). It is recursively presented
  (Theorems~\ref{thm:<>2:pres} and~\ref{thm:=2:pres}), and has
  quadratic growth (Theorem~\ref{thm:<>2:gkdim} and
  Corollary~\ref{cor:=2:gkdim}), so its Gelfand-Kirillov dimension is
  $2$. The ring $\AA$ has an ideal $\KK$, and an embedding
  $\psi:\AA\to M_2(\AA)$, such that all the following:
  $\psi^{-1}:M_2(\KK)\to\KK$, $\KK\to\AA$, $\psi:\AA\to M_2(\AA)$ are
  inclusions with finite cokernel\footnote{i.e.\ the image has finite
    codimension in the target} (Theorem~\ref{thm:galg:jip}).

  Over a field of characteristic $2$, the ring $\AA$ is graded
  (Corollary~\ref{cor:=2:graded}), and may be presented as
  \begin{multline*}
    \AA=\langle A,B,C,D|\,A^2,B^2,C^2,B+C+D,BC,CB,DAD,\\
    \sigma^n(CACACAC),\sigma^n(DACACAD)\text{ for all }n\ge0\rangle,
  \end{multline*}
  where $\sigma$ is the substitution
  $\sigma:\{A,B,C,D\}^*\to\{A,B,C,D\}^*$ defined by
  \[A\mapsto ACA,\quad B\mapsto D,\quad C\mapsto B,\quad D\mapsto C.\]
  The subgroup generated by $\{1+A,1+B,1+C,1+D\}$ is isomorphic to the
  Grigorchuk group $\Gg$. The ring $\AA$ also contains a copy of the
  Laurent polynomials $\F[2][X,X^{-1}]$
  (Theorem~\ref{thm:=2:laurent}).

  If the ground field $\Bbbk$ has characteristic $\neq2$, then $\AA$
  is semiprimitive. If $\Bbbk$ is an algebraic extension of $\F[2]$,
  then $\AA$ is graded nil\footnote{i.e.\ all its homogeneous elements
    are nil}, and its Jacobson radical coincides with its augmentation
  ideal.  If $\Bbbk$ is a non-algebraic extension of $\F[2]$, then
  $\AA$ is a primitive ring, and is not graded nil
  (Theorem~\ref{thm:primjacrad}).
\end{thm}

\noindent The following statement summarizes the main properties of
the rings constructed:
\begin{cor}
  For any field $\Bbbk$, there is a $\Bbbk$-algebra $\KK$ which
  \begin{itemize}
  \item is finitely generated and infinite-dimensional, but has only
    finite-dimensional quotients;
  \item has a subalgebra of finite codimension, isomorphic to
    $M_2(\KK)$;
  \item is prime;
  \item has quadratic growth, and therefore Gelfand-Kirillov dimension
    $2$;
  \item is recursively presented;
  \item satisfies no identity;
  \item contains a transcendental, invertible element;
  \item is semiprimitive if $\Bbbk$ has characteristic $\neq2$;
  \item is graded if $\Bbbk$ has characteristic $2$;
  \item is primitive if $\Bbbk$ is a non-algebraic extension of $\F[2]$;
  \item is graded nil and Jacobson radical if $\Bbbk$ is an algebraic
    extension of $\F[2]$.
  \end{itemize}
\end{cor}

There are interesting examples of primitive, just-infinite algebras
with arbitrary Gelfand-Kirillov dimension~\cite{vishne:primitive};
they are constructed by their presentation (as monomial algebras). The
present construction proceeds in the opposite direction: the algebras
are given as a set of endomorphisms of a vector space, and their
algebraic properties are deduced from the representation.

\subsection{Plan}
Section~\ref{sec:groups} recalls constructions and results concerning
groups acting on rooted trees. A few of the results are new
(Propositions~\ref{prop:gr:hd} and~\ref{prop:gr:growth}); the others
are given with brief proofs, mainly to illustrate the parallelism
between groups and algebras.

Section~\ref{sec:algebras} introduces branch algebras, and develops
general tools and results concerning them; in particular, the branch
algebra associated with a group acting on a rooted tree.

Section~\ref{sec:examples} studies more intricately the branch algebra
associated with the Grigorchuk group. Its study then splits in two
cases, depending on the characteristic being tame ($\neq2$) or wild
($=2$). More results hold in characteristic $2$, in particular because
the branch algebra is graded; some results hold in both cases but the
proofs are simpler in characteristic $2$, and therefore are given in
greater detail there.

\subsection{Notation}
We use the following notational conventions: functions are written
$x\mapsto x^f$ if they are part of a group that acts, and $x\mapsto
f(x)$ otherwise. Generally groups are written in usual capitals ($G$),
and algebras in gothic ($\AA$). We use $\varepsilon$ for the
augmentation map on group rings, $\varpi=\ker\varepsilon$ for the
augmentation ideal, and $\jac\AA$ for the Jacobson radical of $\AA$.

\subsection{Thanks}
We are greatly indebted to Katia Pervova, Said Sidki and Efim Zelmanov
for their open discussions on this topic. Agata Smoktunowicz
generously contributed many interesting remarks concerning the
structure of the Jacobson radical of the rings studied in this paper,
and in particular Lemma~\ref{lem:jac<=>nil}, and Katia Pervova
contributed essential remarks on the nillity of $\AA$. Some of the
results were discovered after experimentation within the computer
algebra system \textsc{Gap}~\cite{gap4:manual}, and its open
development spirit should be commended.  The referee's careful reading
of the paper has been greatly appreciated.

\section{Groups acting on trees}\label{sec:groups}
\subsection{Groups and trees}
We start by reviewing the basic notions associated to groups acting on
rooted trees.
\subsubsection{Trees}
Let $X$ be a set of cardinality $\#X\ge2$, called the \emph{alphabet}.
The \emph{regular rooted tree} on $X$ is $X^*$, the set of (finite)
words over $X$. It admits a natural tree structure by putting an edge
between words of the form $x_1\dots x_n$ and $x_1\dots x_nx_{n+1}$,
for arbitrary $x_i\in X$.  The root is then the empty word.

More pedantically, the tree $X^*$ is the Hasse diagram of the free
monoid $X^*$ on $X$, ordered by right divisibility ($v\le
w\Leftrightarrow\exists u:vu=w$).

Let $G$ be a group with given action on a set $X$. Recall that $A\wr
G$, the \emph{wreath product} of $A$ with $G$, denotes the group
$A^X\rtimes G$, or again pedantically the semi-direct product with $G$
of the sections of the trivial $A$-bundle over $X$.

\subsubsection{Decomposition}\label{ss:gr:decomposition}
Let $W=\aut X^*$ be the group of graph automorphisms of $X^*$. For
each $n\in\N$, the subset $X^n$ of $X^*$ is stable under $W$, and is
called the \emph{$n$th layer} of the tree. The group $W$ admits a
natural map, called the \emph{decomposition}
\[\phi: W\to W\wr\sym X,\]
given by $\phi(g)=(f,\pi_g)$ where $\pi_g\in\sym X$, the
\emph{activity} of $g$, is the restriction of $g$ to the subset
$X\subset X^*$, and $f:X\to W$ is defined by $x^{\pi_g}w^{f(x)}=(xw)^g$,
or in other words $f(x)$ is the compositum $X^*\to
xX^*\stackrel{g}\to\pi_g(x)X^*\to X^*$, where the first and last
arrows are given respectively by insertion and deletion of the first
letter.

The decomposition map can be applied, in turn, to each of the factors
of $W\wr\sym X$. By abuse of notation, we say that we iterate the map
$\phi$ on $W$, yielding $\phi^2:W\to W\wr\sym X\wr\sym X\le
W\wr\sym{X^2}$, etc. More generally, we write $\phi^n:W\to
W\wr\sym{X^n}$, and $\pi^n$ its projection to $\sym{X^n}$.

The action of $W$ on $X^*$ uniquely extends, by continuity, to an
action on $X^\omega$, the (Cantor) set of infinite sequences over $X$.
The self-similarity of $X^\omega$ is expressed via the decomposition
$X^\omega=\bigsqcup_{x\in X}X^\omega$. This gives, for all $n\in\N$, a
continuous map $X^\omega\to X^n$ obtained by truncating a word to its
first $n$ letters.

\subsubsection{$X^*$-bimodule}
There is a left-action $*$ and a right-action $@$ of the free monoid
$X^*$ on $W$, defined for $x\in X$ and $g\in W$ by
\begin{align*}
  x*g&: w\mapsto\begin{cases} x(v^g) & \text{ if }w=xv\\ w & \text{
      otherwise},\end{cases}\\
  g@x&: w\mapsto v\text{ if }(xw)^g=x^gv.
\end{align*}
These actions satisfy the following properties:
\begin{xalignat}{2}
  (g@v)@w &= g@(vw), &       v*(w*g) &= (vw)*g,\\
  (gh)@v  &= (g@v)(h@v^g), & v*(gh)  &= (v*g)(v*h),\\
  g       &= (v*g)@v, &      g       &= \left(\prod_{v\in X^n}v*(g@v)\right)\,\pi^n_g,
\end{xalignat}
where in the last expression the $v*(g@v)$ mutually commute when $v$
ranges over the $n$th layer $X^n$.

In this terminology, when we wrote the decomposition as $\phi(g) =
(f,\pi_g)$, we had $f(x)=g@x$.

\subsubsection{Branchness}
Let $G<W$ be a group acting on the regular rooted tree $X^*$. The
\emph{vertex stabilizer} $\stab_G(v)$ is the subgroup of
$G$ fixing $v\in X^*$. The group $G$ is
\begin{description}
  \def\makelabel#1{\textbf{\emph{#1}},}
\item[level-transitive] if $G$ acts transitively on $X^n$ for all
  $n\in\N$;
\item[recursive] if $G@x<G$ for all $x\in X$;
\item[weakly recurrent] if it is level-transitive, and $G@x=G$ for all
  $x\in X$;
\item[recurrent] if it is level-transitive, and $\stab_G(x)@x=G$ for
  all $x\in X$;
\item[weakly branch] if $G$ is level-transitive, and $(v*G)\cap G$ is
  non-trivial for all $v\in X^*$;
\item[weakly regular branch] if $G$ is level-transitive, and has a
  non-trivial normal subgroup $K$, called the \emph{branching
    subgroup}, with $x*K<K$ for all $x\in X$;
\item[branch] if $G$ is level-transitive, and $\langle (v*G)\cap
  G:v\in X^n\rangle$ has finite index in $G$ for all $n\in\N$;
\item[regular branch] if $G$ is level-transitive, and has a
  finite-index normal subgroup $K$ with $x*K<K$ for all $x\in X$;
\end{description}

Weak branchness can be reformulated in terms of the action on
$X^\omega$. Then $G$ is weakly branch if every closed set
$F\subsetneqq X^\omega$ has a non-trivial fixator
$\operatorname{Fix}_G(F)=\setsuch{g\in G}{g(f)=f\,\forall f\in F}$.

Remark that if $G$ is branch, then $K^X$ has finite index in
$\phi(K)$, because it has finite index in $G^X$ and in $G\wr\sym X$.

Remark also that if $G$ is weakly regular branch, then there is a
unique maximal branching subgroup $K$; it is
\[K=\bigcap_{v\in X^*}(G\cap(v*G))@v.\]

\begin{prop}\label{prop:xtive}
  If $G$ is transitive on $X$ and $\stab_G(x)@x<G$ for all $x\in X$,
  then it is recurrent, and therefore its action on
  $(X^\omega,\text{Bernoulli})$ is ergodic.  In particular, $G$ is
  infinite.
\end{prop}
\begin{proof}
  Proceed by induction on $n$.  Consider a layer $X^n$ of the tree, and two
  vertices $x_1\dots x_n$ and $y_1\dots y_n$. Since $G$ is branch, it acts
  transitively on $X$, so $x_1\dots x_n$ and $y_1x_2\dots x_n$ belong to the
  same orbit. By induction, $x_2\dots x_n$ and $y_2\dots y_n$ are in the same
  $G$-orbit; therefore, since $\stab_G(y_1)=G$, the vertices $y_1x_2\dots x_n$
  and $y_1\dots y_n$ belong to the same orbit.
  
  If the action is not ergodic, let $A\subset X^\omega$ be an
  invariant subset of non-$\{0,1\}$ measure. Then there exists
  $n\in\N$ such that $X^\omega\to X^n$ is not onto; its image is a
  $G$-orbit, and thus the action of $G$ is not transitive on the $n$th
  layer.
\end{proof}

%

\begin{prop}
  If $G$ is regular branch, then it is regular weakly branch and
  branch; if it is branch, then it is weakly branch; if it is regular
  weakly branch, then it is weakly branch.
\end{prop}
\begin{proof}
  Let $G$ be a regular branch group, with branching subgroup $K$. By
  Proposition~\ref{prop:xtive}, $G$ is infinite so $K$ is non-trivial.
  This shows that $G$ is regular weakly branch.  Assume now only that
  $K$ is non-trivial, and let $v\in X^n$ be any vertex.  Since
  $K^{X^n}\leq\phi^n(G)$, we may take any $k\neq1$ in $K$ and consider
  the element $k*v\in G$. This shows that $G$ is weakly branch. The
  other implications are of the same nature.
\end{proof}

Note finally that the group $G$ is determined by a generating set $S$
and the restriction of the decomposition map $\phi$ to $S$, in the
following sense:
\begin{prop}\label{prop:unique}
  Let $F$ be a group generated by a set $S$, and let $\phi:F\to
  F\wr\sym X$ be any map. Then there exists a unique subgroup $G$ of
  $W=\aut X^*$ that is generated by $S$ and has decomposition map
  induced by $\phi$ through the canonical map $F\to G$.
\end{prop}
\begin{proof}
  The decomposition map $\phi$ yields, by iteration, a map
  $F\to\sym{X^n}$ for all $n\in\N$. This defines an action on the
  $n$th layer of the tree $X^*$, and since they are compatible with
  each other they define an action of $F$ on the tree. We let $G$ be
  the quotient of $F$ by the kernel of this action. On the other hand,
  the action of the generators, and therefore of $G$, is determined by
  $\phi$, so $G$ is unique.
\end{proof}
In particular, $F$ may be the free group on $S$, and $\phi$ may be
simply defined by the choice, for each generator in $S$, of $\#X$
words and a permutation.

Therefore, in defining a recurrent group, we will only give a list of
generators, and their images under $\phi$. If $X=\{1,\dots,q\}$, we
describe $\phi$ on generators with the notation
\[\phi(g)=\pair<g@1,\dots,g@q>\pi_g,\text{ or even
}\phi(g)=\pair<g_1,\dots,g_q>\text{ if }\pi_g=1,
\]
rather than in the form $\phi(g)=(f,\pi)$ with $f(x)=g@x$.

Note that there may exist other groups $G'$ generated by $S$, and such
that the natural map $F\to G'$ induces an injective map $G'\to
G'\wr\sym X$. However, such $G'$ will not act \emph{faithfully} on
$X^*$. The group $G$ defined by Proposition~\ref{prop:unique} is the
\emph{smallest} quotient of $F$ through which the decomposition map
factors.

Weakly branch groups $G$ are known to satisfy no identity; i.e.\ for
every $w\neq1$ in the free group $F(y_1,\dots,y_k)$ there exist
$g_1,\dots,g_k\in G$ with $w(g_1,\dots,g_k)\neq1$. We quote the
following general result, due to Mikl\'os Ab\'ert:
\begin{prop}[\cite{abert:nonfree}*{Theorem~1}]\label{prop:abert}
  Let $G$ be a group acting on a set $X$, such that for every finite $Y\subset
  X$ the fixator\footnote{aka ``pointwise stabilizer''} of $Y$ does not fix any
  point in $X\setminus Y$. Then $G$ does not satisfy any identity.
\end{prop}
His proof goes as follows: let $w_i$ be the length-$i$ prefix of $w$,
and let $x\in X$ be any. Then, inductively on $i$, one shows that
there exist $g=(g_1,\dots,g_k)\in G^k$ such that
$x,x^{w_1(g)},\dots,x^{w_i(g)}$ are all distinct. The following is a
weakening of~\cite{abert:nonfree}*{Corollary~4}.

\begin{cor}\label{cor:nolaw}
  If $G$ is weakly branch, then it does not satisfy any identity.
\end{cor}
\begin{proof}
  Let $G$ act on the boundary $X^\omega$ of the tree $X^*$. Let $Y\subset X$ be
  a finite subset, and let $\xi\in X\setminus Y$ be any. Then there exists a
  vertex $v\in X^*$ on the geodesic $\xi$ but on none of the geodesics in $Y$.
  Set $K=G\cap(G*v)$. Since $G$ is weakly branch, $K$ is non-trivial. Assume by
  contradiction that $K$ fixes $\xi$. Then since $K$ is invariant under the
  stabilizer of $v$, and $G$ acts level-transitively, it follows that $K$ also
  fixes all images of $\xi$ under the stabilizer of $\xi$; this is a dense
  subset of $vX^\omega$, so $K$ fixes $X^\omega$, which contradicts the
  non-triviality of $K$. Therefore there exists $g\in K$ with $g|_Y=1$ and
  $\xi^g\neq\xi$, so the conditions of Proposition~\ref{prop:abert} are
  satisfied.
\end{proof}

\begin{prop}[\cite{lavreniuk-n:rigidity}*{Lemma~5.4}]\label{prop:wb:centerless}
  Let $G$ be a weakly branch group. Then its centre is trivial.
\end{prop}
\begin{proof}
  More generally, take $g\neq1\in\aut(X^*)$; then it moves a vertex
  $u$. Since $G$ is weakly branch, there is $h\neq1$ acting only on
  the subtree $uX^*$, and $[g,h]\neq1$.
\end{proof}

\subsection{Dimension}
Every countable residually-$p$ group has a representation as a
subgroup of $\aut X^*$, for $X=\{1,\dots,p\}$: fix a descending
filtration $G=G_0\ge G_1\ge G_2\ge\dots$ with $\bigcap G_n=\{1\}$ and
$[G_n:G_{n+1}]=p$; identify $X$ with $G_n/G_{n+1}$. Then $G/G_n$ is
identified with $X^n$, and $G$ acts faithfully, by multiplication on
cosets, on the tree $X^*$. In general, this action will not be
recurrent. Moreover, this action may be ``inefficient'' in that the
quotient of $G$ represented by the action on $X^n$ may be quite small
--- if $G_n\triangleleft G$ this quotient is $G/G_n$ of order $p^n$,
while the largest $p$-group acting on $X^n$ has order
$p^{(p^n-1)/(p-1)}$. This motivates the following definition.

Let $W_n=\pi^n(W)$ be quotient of $W$ acting on $X^n$. We give $W$ the
structure of a compact, totally disconnected metric space by setting
\[d(g,h)=\inf\setsuch{1/\#W_n}{\pi^n(g)=\pi^n(h)}.\] We obtain in this
way the notion of \emph{closure} and \emph{Hausdorff dimension}.
Explicitly, for a subgroup $G<W$, we have
by~\cite{abercrombie:subgroups}
\[\Hdim(G)=\liminf_{n\to\infty}\frac{\log\#\pi^n(G)}{\log\#W_n};\]
see also~\cite{barnea-s:hausdorff}. The Hausdorff dimension of $G$
coincides with that of its closure.

\subsubsection{The tree closure}\label{ss:gr:tc}
Let $P\le\sym X$ be any group acting on $X$. The \emph{tree closure}
of $P$ is the subgroup $\overline P$ of $W$ consisting of all $g\in W$
such that $\pi^n(g)\in\wr^n P\le\sym{X^n}$ for all $n\in\N$. It is the
inverse limit of the groups $\wr^n P$, and is a closed subgroup of
$W$.

We have $\overline P=\overline P\wr P$, and $\pi^n(\overline
P)=\pi^{n-1}(\overline P)\wr P$, so $\#\pi^n(\overline
P)=(\#\pi^{n-1}(\overline P))^{\#X}\#P$, and therefore
\begin{equation}\label{eq:neumanndim}
  \#\pi^n(\overline P)=(\#P)^{\frac{\#X^n-1}{\#X-1}}.
\end{equation}
In particular, $\#W_n=(\#\sym X)^{(\#X^n-1)/(\#X-1)}$, and $\overline
P$ has Hausdorff dimension $\log\#P/\log(\#X!)$.

If $p$ is prime and $X=\{1,\dots,p\}$, we will often consider
subgroups $G$ of $W_p=\overline P$, where
$P=\langle(1,2,\dots,p)\rangle$ is a $p$-Sylow of $\sym X$. The
dimension of $G$ will be then computed relative to $W_p$, by the
simple formula
\[\Hdim_p(G)=\frac{\Hdim G}{\Hdim W_p}=\Hdim G\frac{\log(p!)}{\log p}.\]

\begin{prop}\label{prop:gr:hd}
  Let $G$ be a regular branch group. Then $G$ has positive Hausdorff
  dimension.
  
  If furthermore $G$ is a subgroup of $W_p$, then its relative
  Hausdorff dimension $\Hdim_p$ is rational.
\end{prop}
\begin{proof}
  Let $G$ have branching subgroup $K$, and for all $n\in\N$ set
  $G_n=\pi^n(G)$. Let $M\in\N$ be large enough so that
  $G/\phi^{-2}(K^{X^2})$ maps isomorphically into
  $G_M/\pi^{M-2}(K^{X^2})$. We then have, for all $n\ge M$,
  \begin{align}
    \#G_n&=[G:K]\#\pi^n(K)=[G:K][\phi(K):K^X](\#\pi^{n-1}(K))^{\#X}\label{eq:poshd}\\
    &=[G:K]^{1-\#X}[\phi(K):K^X](\#G_{n-1})^{\#X}.\notag
  \end{align}
  Write $\log\#G_n=\alpha\#X^n+\beta$, for some $\alpha,\beta$ to
  be determined; we have, again for $n\ge M$,
  \[\alpha\#X^n+\beta=(1-\#X)\log[G:K]+\log[\phi(K):K^X]+\#X(\alpha\#X^{n-1}+\beta),\]
  so $\beta=\log[G:K]-\log[\phi(K):K^X]/(\#X-1)$. Then set
  $\alpha=(\log\#G_M-\beta)/\#X^M$. We have solved the recurrence for
  $\#G_n$, and $\alpha>0$ because $G_n$ has unbounded order.

  Now it suffices to note that $\Hdim(G)=\alpha(\#X-1)/\log(\#X!)$ to
  obtain $\Hdim(G)>0$.
  
  For the last claim, note that all indices in~\eqref{eq:poshd} are
  powers of $p$, and hence their logarithms in base $p$ are integers.
\end{proof}

\begin{question}
  Mikl\'os Ab\'ert and B\'alint Vir\'ag~\cite{abert-v:dimension} show
  that there exist free subgroups of $W$ of Hausdorff dimension $1$.
  Is there a finitely generated recurrent group of dimension $1$? A
  branch group?
\end{question}

\subsection{Growth}
Let $G$ be a group generated by a finite set $S$. The \emph{length} of
$g\in G$ is defined as $\|g\|=\min\setsuch{n}{g=s_1\dots s_n\text{ for
    some }s_i\in S}$. The \emph{word growth} of $G$ is the function
\[f_{G,S}(n)=\#\setsuch{g\in G}{\|g\|\le n}.\]
This function depends on the choice of generating set $S$. Given
$f,g:\N\to\R$, say $f\precsim g$ if there exists $M\in\N$ with
$f(n)\le g(Mn)$, and say $f\sim g$ if $f\precsim g\precsim f$; then
the equivalence class of $f_{G,S}$ is independent of $S$. The group
$G$ has \emph{exponential growth} if $f_{G,S}\sim e^n$, and polynomial
growth if $f_{G,S}\precsim n^D$ for some $D\in\N$. In all other cases,
$f_{G,S}$ grows faster than any polynomial and slower than any
exponential, and $G$ has \emph{intermediate growth}.  If furthermore
$f_{G,S}(n)\ge A^n$ for some $A>1$, uniformly on $S$, then $G$ has
\emph{uniformly exponential growth}.

More generally, let $E$ be a space on which $G$ acts, and let $*\in E$
be any. Then the \emph{growth of $E$} is the function
\[f_{E,*,S}(n)=\#\setsuch{e\in E}{e=g*\text{ with }\|g\|\le n}.\]
If $E=G$ with left regular action, we recover the previous definition
of growth. We will be interested in the case $E=X^\omega$ with the
natural action of $G$, or equivalently of $E=G/\stab_G(*)$ for some
$*\in X^\omega$.

\subsubsection{Contraction}\label{sss:contraction}
Let $G$ be a finitely generated recurrent group. It is
\emph{contracting} if there exist $\lambda<1$, $n\in\N$ and $K$ such
that, for all $g\in G$ and $v\in X^n$ we have
$\|g@v\|\le\lambda\|g\|+K$.

\begin{prop}[\cite{bartholdi-g-n:fractal}, Proposition~8.11]\label{prop:contr}
  If $G$ is contracting, then the growth of $(X^\omega,*)$ is
  polynomial, of degree at most $-n\log\#X/\log\lambda$.
  
  Conversely, if $(X^\omega,*)$ has polynomial growth of degree $d$,
  then $G$ is contracting for any $n$ large enough and any
  $\lambda>(\#X)^{-n/d}$.
\end{prop}

\begin{prop}\label{prop:gr:growth}
  If $G$ is a finitely generated branch group, and $(X^\omega,*)$ has
  polynomial growth of degree $d$, then $G$ has growth
  \[f_G(n)\succsim\exp\left(n^{d/(d+1)}\right).\]
\end{prop}
\begin{proof}
  Let us write $q=\#X$. Let $K$ be a branching subgroup, and set
  $R_0=\min\setsuch{\|g\|}{g\in K,g\neq1}$. Let $n\in\N$ and $v\in
  X^n$ be given. For $g\in K$ satisfying $\|g\|\le R_0$, set
  $h_{v,g}=v*g$. By Proposition~\ref{prop:contr}, we have
  $\|h_{v,g}\|\le q^{n/d}\|g\|\le q^{n/d}R_0$.
  
  We now choose for all $v\in X^n$ some $g_v\in K$ with $\|g_v\|\le
  R_0$, and consider the corresponding element
  \[h=\prod_{v\in X^n}h_{v,g_v}.\]
  On the one hand, there are at least $2^{q^n}$ such elements, because there
  are at least $2$ choices for each $g_v$. On the other hand, such an element
  has length at most $q^nq^{n/d}R_0$. If $f(R)$ denote the growth function of
  $G$, we therefore have $f(q^{n+n/d}R_0)\ge 2^{q^n}$, or in other words
  \[f(R)\succsim\exp\left(q^{\log R/(1+\frac1d)\log q}\right)=\exp\left(R^{d/(d+1)}\right).\]
\end{proof}

\subsection{Main examples}\label{ss:gr:ex}
The first example of a branch group is $W$ itself, with branching
subgroup $K=W$. In this paper, however, we are mainly concerned with
countable groups. Assume therefore that $X$ is finite, and choose a
section $\widetilde{\sym X}$ of $\pi:W\to\sym X$, for instance lifting
$\rho\in\sym X$ to $\widetilde\rho:x_1\dots
x_n\mapsto\rho(x_1)x_2\dots x_n$. The \emph{finitary group} $\sym
X^\flat$ is the subgroup of $W$ generated by the
$\phi^{-n}\big(\widetilde{\sym X}^{X^n}\big)$, for all $n\in\N$. It is
locally finite.

More generally, let $P$ be the lift to $W$ of a transitive subgroup of
$\sym X$. The \emph{finitary closure} of $P$ is then the subgroup
$P^\flat$ of $W$ generated by the $\phi^{-n}(P^{X^n})$, for all
$n\in\N$. If $P$ is countable, then $P^\flat$ is a countable subgroup
of the tree closure $\overline P$ of $P$.

Much of the interest in branch groups comes from the fact that
finitely generated examples exist. The most important ones are:
\subsubsection{The Neumann groups}
Take $P$ a perfect, $2$-transitive subgroup of $\sym X$, and choose
$a,b\in X$. Consider two copies $P$, $\overline P$ of $P$, and let
them act on $X^*$ as follows:
\[(x_1\dots x_n)^p=(x_1^p)x_2\dots x_n;\quad
(x_1\dots x_n)^{\overline p}=\begin{cases}
  x_1(x_2\dots x_n)^p & \text{ if }x_1=a,\\
  x_1(x_2\dots x_n)^{\overline p} & \text{ if }x_1=b,\\
  x_1\dots x_n & \text{ else}.
\end{cases}
\]
Let $G$ be the group generated by these two images of $P$. Then $G$ is
a perfect group, studied by Peter Neumann in~\cite{neumann:pride}; it
is branch, with branching subgroup $K=G$. Indeed choose $r,s\in P$
with $a^r=a\neq a^s$ and $b^r\neq b=b^s$. Then $\phi[\overline
P,\overline P^r]=P\times1\times\dots\times1$ and $\phi[\overline
P,\overline P^s]=1\times\overline P\times\dots\times1$, so
$\phi(\overline P^P)$ contains $G\times1\dots\times1$ and therefore
contains $G\times\dots\times G$. Note that $P$ is isomorphic to
$\phi(G)/G^X$.
  
The group $G$ is more simply defined by its decomposition map: $G$ is
the unique subgroup of $W$ generated by two copies $P\sqcup\overline
P$ of $P$ and with decompositions
\[\phi(p)=\pair<1,\dots,1>p,\quad
\phi(\overline p)=\pair<p,\overline p,1,\dots,1>,
\]
with in the last expression the `$p$' in position $a$ and the
`$\overline p$' in position $b$.

The example $P=\mathrm{PSL}_3(2)$, in its action on the $7$-point
projective plane, was considered in~\cite{bartholdi:nueg}, where $G$
was shown to have non-uniformly exponential word growth; see
also~\cite{wilson:ueg}. These groups are contracting with $n=1$ and
$\lambda=\frac12$.
  
The Hausdorff dimension of $G$ is $\log\#P/\log(\#X!)$,
by~\eqref{eq:neumanndim}.

\subsubsection{The Grigorchuk group}\label{sss:gr:gg}
This group $\Gg$ acts on the binary tree, with $X=\{1,2\}$. It is best
described as the group generated by $\{a,b,c,d\}$, with given
decompositions:
\begin{equation}\label{eq:grgr:decomp}
  \phi(a)=\pair<1,1>(1,2),\quad\phi(b)=\pair<a,c>,\quad\phi(c)=\pair<a,d>,\quad\phi(d)=\pair<1,b>.
\end{equation}
This group was studied by Rostislav Grigorchuk, who showed
in~\cite{grigorchuk:burnside} that $\Gg$ is a f.g.\ infinite torsion
group --- also known as a ``Burnside group''. He then showed
in~\cite{grigorchuk:growth} that it has word-growth intermediate
between polynomial and exponential; the more precise bounds
\[\exp(n^{0.5157})\precsim f_G\precsim\exp(n^{0.7675})\]
appear in~\cites{bartholdi:upperbd,bartholdi:lowerbd}. This group is
contracting with $n=1$ and $\lambda=\frac12$.

The group $\Gg$ is a branch group, with branching subgroup
$K=\langle[a,b]\rangle^\Gg$ of index $16$. Indeed set $x=[a,b]$; then
$\phi[x^{-1},d]=\pair<1,x>$ so $\phi(K)$ contains $K\times K$. Set
$x=[a,b]$; then, as a group, $K$ is generated by
$\{x,[x,d],[x,d^a]\}$.

The finite quotient $\pi^n(\Gg)$ has order $2^{5\cdot 2^{n-3}+2}$, for
$n\ge3$. It follows that $\Gg$ has Hausdorff dimension $5/8$.
  
Igor Lys\"enok obtained in~\cite{lysionok:pres} a presentation of
$\Gg$ by generators and relations:
\begin{prop}[\cite{lysionok:pres}]\label{prop:lysenok}
  Consider the endomorphism $\sigma$ of $\{a,b,c,d\}^*$ defined by
  \begin{equation}\label{eq:sigma}
    a\mapsto aca,\quad b\mapsto d,\quad d\mapsto c,\quad c\mapsto b.
  \end{equation}
  Then
  \begin{equation}
    \Gg=\big\langle a,b,c,d\big|a^2,b^2,c^2,d^2,bcd,
    \sigma^n(ad)^4,\sigma^n(adacac)^4\;\forall n\ge0\big\rangle.\label{eq:lysenok}
  \end{equation}
\end{prop}
Note that if the relator $r$ is understood as $r=1$, this gives a ring
presentation of the group ring $\Bbbk\Gg$. Since the algebra $\AA$
mentioned in the introduction is a quotient of $\Bbbk\Gg$, it must
have stronger relations than the above --- see
Theorem~\ref{thm:<>2:pres}.  The last two families of relations, in
$\Bbbk\Gg$, may be written as $\sigma^n[d^a-1,d-1]=0$ and
$\sigma^n[d^{(ac)^2a}-1,d-1]=0$. In essence, these relations are
strengthened in $\AA$ to $\sigma^n((d^a-1)(d-1))=0$ and
$\sigma^n((d^{(ac)^2a}-1))(d-1)=0$ respectively.

\subsubsection{The Gupta-Sidki group}\label{sss:gs}
This group $\GS$ acts on the ternary tree, with $X=\{1,2,3\}$. It is
best described as the group generated by $\{x,\gamma\}$, with
decompositions
\[\phi(x)=\pair<1,1,1>(1,2,3),\quad \phi(\gamma)=\pair<\gamma,x,x^{-1}>.\]
This group was studied by Narain Gupta and Said
Sidki~\cite{gupta-s:burnside}, who showed that $\GS$ is an infinite
$3$-torsion group.
  
This group is contracting with $n=1$ and $\lambda=\frac12$.
  
The finite quotient $\pi^n(\GS)$ has order $3^{2\cdot 3^{n-1}+1}$, for
$n\ge2$. It follows that $\GS$ has Hausdorff dimension $4/9$ in $W_3$.

The group $\GS$ is branch, with branching subgroup $\GS'=[\GS,\GS]$.
Indeed $\phi(\GS')$ contains $\GS'\times\GS'\times\GS'$, because
$\phi([\gamma^{-1}\gamma^{-x^2},\gamma^x\gamma])=\pair<1,1,[x,\gamma]>$.
  
Later Said Sidki constructed a presentation of $\GS$ by generators and
relations~\cite{sidki:pres}, and associated an algebra to $\GS$ ---
see Theorem~\ref{thm:sidki}.

\subsubsection{Weakly branch groups}
Most known examples of recurrent groups are weakly branch. Among those
that are not branch, one of the first to be considered acts on the
ternary tree $\{1,2,3\}^*$:
\[\BG=\langle x,\delta\rangle\text{ given by }\phi(x)=\pair<1,1,1>(1,2,3),\quad\phi(\delta)=\pair<\delta,x,x>;\]
It was studied along with $\Gg$, $\GS$ and two other examples
in~\cites{bartholdi-g:parabolic,bartholdi-g:spectrum}.  The finite
quotient $\pi^n(\BG)$ has order $3^{\frac14(3^n+2n+3)}$, for $n\ge2$.
It follows that $\BG$ has Hausdorff dimension $1/2$ in $W_3$.
  
Two interesting examples, acting on the binary tree, were also found:
\subsubsection{The ``BSV'' group}
\[G_1=\langle \tau,\mu\rangle\text{ given by
}\phi(\tau)=\pair<1,\tau>(1,2),\quad\phi(\mu)=\pair<1,\mu^{-1}>(1,2);
\]
it was studied in~\cite{brunner-s-v:nonsolvable}, who showed that it
is torsion-free, weakly branch, and constructed a presentation of
$G_1$.  The finite quotient $\pi^{2n}(G_1)$ has order
$2^{\frac13(2^{2n}-1)+n}$, for $n\ge1$.  It follows that $G_1$ has
Hausdorff dimension $1/3$.

\subsubsection{The Basilica group}
\[G_2=\langle a,b\rangle\text{ given by
}\phi(a)=\pair<1,b>(1,2),\quad\phi(b)=\pair<1,a>;
\]
it was studied in~\cite{grigorchuk-z:torsionfree}, who showed that it
is torsion-free and weakly branch, and
in~\cite{bartholdi-v:amenability}, who showed that it is amenable,
though not ``subexponentially elementary amenable''. The finite
quotient $\pi^{2n}(G_2)$ has order $2^{\frac23(2^{2n}-1)+n}$, for
$n\ge1$.  It follows that $G_2$ has Hausdorff dimension $2/3$.
  
All of these groups are contracting with $n=1$ and $\lambda=\frac1{\sqrt2}$.

\subsubsection{The odometer}\label{sss:odo}
This is a group acting on $\{1,2\}^*$:
\[\Z=\langle\tau\rangle,\quad\phi(\tau)=\pair<1,\tau>(1,2).\]
Its action on the $n$th layer is via a $2^n$-cycle. It is not weakly
branch.
\subsubsection{The Lamplighter group}\label{sss:ll}
This is the group $G=(\Z/2)^{(\Z)}\rtimes\Z$, the semidirect product
with $\Z$ of finitely-supported $\Z/2$-valued functions on $\Z$. It
acts on $\{1,2\}^*$:
\[G=\langle a,b\rangle,\quad\phi(a)=\pair<a,b>(1,2),\quad\phi(b)=\pair<a,b>.\]
Again this group is not weakly branch.

%

\section{Algebras}\label{sec:algebras}
We consider various definitions of ``recurrence'' and ``branchness''
in the context of algebras. Let $\Bbbk$ be a field, fixed throughout
this section.
\subsection{Associative algebras}
If $X$ is a set, we write $M_X(\Bbbk)=M_X$ the matrix algebra of
endomorphisms of the vector space $\Bbbk X$, and for a $\Bbbk$-algebra
$\AA$ we write $M_X(\AA)=M_X(\Bbbk)\otimes\AA$.

\subsubsection{Recurrent transitive algebras}\label{sss:rt}
A \emph{recurrent transitive algebra} is an associative algebra $\AA$,
given with an injective homomorphism $\psi:\AA\to M_X(\AA)$, for some
set $X$, such that for every $x,y\in X$ the linear map $\AA\to
M_X(\AA)\to \AA$, obtained by projecting $\psi(\AA)$ on its $(x,y)$
matrix entry, is onto.

The map $\psi$ is called the \emph{decomposition} of $\AA$, and can be
iterated, yielding a map $\psi^n:\AA\to M_{X^n}(\AA)$.

The most naive examples are as follows: consider the vector space
$V=\Bbbk X^\omega$, and $\AA=\End(V)$. The decomposition map is given by
$\psi:a\mapsto(a_{x,y})$ where $a_{x,y}$ is defined on the basis vectors
$w\in X^*$ as follows: if $a(xw)=\sum b_vv$, then
\[a_{x,y}(w)=\sum_{v=yv'\in X^\omega} b_vv'.\]

Similarly, consider the vector space $V=\Bbbk^{X^\omega}$ of functions
on $X^\omega$, and $\AA=\End(V)$. The decomposition map is given by
\[\psi(a)=(a_{x,y})\text{ where }a_{x,y}(f)(w)=a(v\mapsto f(xv))(yw).\]

These examples are meant to illustrate the connection between
action on $X^\omega$ and recurrent algebras; they will not be
considered below. However, all our algebras will be subalgebras of
these, i.e.\ contained in $\otimes^\omega M_X=M_{X^\omega}$.

\subsubsection{Decomposition}
Similarly to Proposition~\ref{prop:unique}, a recurrent transitive
algebra may be defined by its decomposition map, in the following
sense:
\begin{lem}
  If $\FF$ is an algebra generated by a set $S$, and $\psi:\FF\to M_X(\FF)$
  is a map such that $\psi(\FF)_{x,y}=\FF$ for all $x,y\in X$, then there
  exists a unique minimal quotient of $\FF$ that is a recurrent
  transitive algebra.
\end{lem}
\begin{proof}
  Set $\II_0=\ker\psi$ and $\II_{n+1}=\psi^{-1}M_X(\II_n)$ for
  $n\in\N$ and $\II=\bigcup_{n\in\N}\II_n$. Then $\II$ is an ideal in
  $\FF$, and $\FF/\II$ is a recurrent transitive algebra. Consider the
  ideal $\JJ$ generated by all ideals $\KK\le\FF$ such that
  $\psi(\KK)\le M_X(\KK)$; then $\II\le\JJ$, and $\AA=\FF/\JJ$ is the
  required minimal quotient of $\FF$.
\end{proof}
It follows that a branch algebra may be defined by a choice, for each
generator in $S$, of $\#X^2$ elements of the free algebra
$\Bbbk\langle S\rangle$.  Note that we do not mention any topology on
$\AA$; if $\AA$ is to be, say, in the category of $C^*$-algebras, then the
definition becomes much more intricate due to the absence of free
objects in that category. The best approach is probably that of a
$C^*$-bimodule considered in~\cite{nekrashevych:bimodule}.

An important feature is missing from the algebras of
\S\ref{sss:rt}, namely the existence of finite-dimensional
quotients similar to group actions on layers. These are introduced as
follows:
\subsubsection{Augmented algebras}
Let $\AA$ be a recurrent transitive algebra. It is \emph{augmented} if
there exists a homomorphism $\varepsilon:\AA\to\Bbbk$, called the
\emph{augmentation}, and a subalgebra $\PP$ of $M_X$ with a
homomorphism $\zeta:\PP\to\Bbbk$, such that the diagram
\[\xymatrix{{\AA}\ar[r]^{\psi}\ar[d]_{\varepsilon} &
  {\psi(\AA)\makebox[0mm][l]{${}\le M_X\otimes \AA$}}\ar[d]^{1\otimes\varepsilon}\\
  {\Bbbk} & {\PP}\ar[l]^{\zeta}}
\]
commutes.  We abbreviate ``augmented recurrent transitive algebra'' to
\emph{art algebra}, or \emph{$\PP$-art algebra} if we wish to
emphasize which $\PP\le M_X$ is used.

Let $\PP$ be a subalgebra of $M_X$, with augmentation
$\zeta:\PP\to\Bbbk$.  There are two fundamental examples of art
algebras, constructed as follows:

\subsubsection{The ``tree closure'' $\overline \PP$}\label{sss:treeclosure}
We define for all $n\in\N$ an augmented algebra $\PP_n\le M_{X^n}$,
with $\zeta_n:\PP_n\to\Bbbk$, for $n\in\N$ by $\PP_1=\PP$,
$\zeta_1=\zeta$, and
\[\PP_{n+1}=\langle m\otimes p\in M_X\otimes \PP_n|\,\zeta_n(p)m\in \PP\rangle.\]
Its augmentation is given by $\zeta_{n+1}(m\otimes
p)=\zeta(\zeta_n(p)m)$.

Then there is a natural map $\PP_{n+1}\to \PP_n$, defined by
$m_1\otimes\dots\otimes m_{n+1}\mapsto
\zeta(m_{n+1})m_1\otimes\dots\otimes m_n$. We set
$\overline\PP=\varprojlim \PP_n$.
  
Then $\overline \PP$ is an art algebra: for $a\in\overline \PP$, write
$a=\varprojlim a_n$ with $a_n\in \PP_n$. Then $a_{n+1}=\sum m_n\otimes
p_n$ with $m_n\in M_X$ and $p_n\in \PP_n$. The sequence $m_n$ is
constant equal to $m$, and we set $\psi(a)=\sum m\otimes\varprojlim
p_n$.

The following diagram gives a natural map $\AA\to\overline\PP$ for any
$\PP$-art algebra $\AA$. We will always suppose that this map is injective.
\[\xymatrix{{\AA}\ar[r]^{\psi}\ar[d]_{\varepsilon} &
  {\psi\AA}\ar[r]^{\psi}\ar[d]^{1\otimes\varepsilon} &
  {\psi^2\AA}\ar[r]^{\psi}\ar[d]^{1\otimes1\otimes\varepsilon} & {\dots}\\
  {\Bbbk} & {\PP}\ar[l]_{\zeta} & {\PP_2}\ar[l]_{1\otimes\zeta} &
  {\dots}\ar[l] & {\overline\PP}\ar[l]}\]
  
\subsubsection{The ``finitary closure''}
This construction starts as above, by noting that the map
$\PP_{n+1}\to \PP_n$, $a_n\otimes p\mapsto \zeta(p)a_n$, is split by
$a_n\mapsto a_n\otimes 1$. We let $\PP^\flat$ be the direct limit of
the $\PP_n$'s along these inclusions.

Then $\PP^\flat$ is also an art algebra. Its decomposition is defined
on $\PP_n$ as above: $\psi(m\otimes p)=m\otimes p$ for $m\in M_X,p\in
\PP_n,m\otimes p\in \PP_{n+1}$.

In some sense, $\overline \PP$ is the maximal $\PP$-art algebra, and
$\PP^\flat$ is a minimal $\PP$-art algebra. More precisely:
\begin{prop}\label{prop:alg:unique}
  Let $\FF$ be an augmented algebra generated by a set $S$, and let
  $\psi:\FF\to M_X(\FF)$ be a map such that $\psi(\FF)_{x,y}=\FF$ for
  all $x,y\in X$. Set $\PP=\varepsilon\psi(\FF)\le M_X$, and
  assume that the augmentation $\epsilon:\FF\to\Bbbk$ factors to
  $\zeta:\PP\to\Bbbk$.

  Then there exists a unique art subalgebra $\AA$ of $\overline \PP$ that
  is generated by $S$ and has decomposition map induced by $\psi$
  through the canonical map $\FF\to \AA$.
\end{prop}
\begin{proof}
  For all $n\in\N$ there exists a map
  $\pi^n=\varepsilon\psi^n:\FF\to\PP_n$, and these maps are compatible
  in that $(1^{\otimes n}\otimes\zeta)\pi^{n+1}=\pi_n$. There is
  therefore a map $\pi:\FF\to\overline\PP$, and we let $\AA$ be the
  image of $\pi$. This proves the existence part.

  Let $\AA'=\FF/\JJ'$ be another image of $\FF$ in $\overline\PP$.
  Write $\JJ=\ker\pi$. Then by definition of art algebra the images of
  $\AA$ in $\PP_n$ must be $\pi^n(\FF)$, so $\JJ'\le\ker\pi^n$,
  and $\JJ'\le\JJ$. It follows that $\JJ'=\JJ$, because $\AA$ and
  $\AA'$ are both contained in $\overline\PP$.
\end{proof}

If $X=\{1,\dots,q\}$, then a maximal augmented subalgebra of $M_X$ is
$\PP\cong M_{q-1}\oplus\Bbbk$, where the augmentation vanishes on
$M_{q-1}$. The examples we shall consider fall into this class.

For $V$ a vector space, we denote by $V^\circ$ its dual, and we
consider $V\otimes V^\circ$ as a subspace of $\End(V)$, under the
natural identification $(v\otimes\xi)(w)=\xi(w)\cdot v$.

\subsubsection{Branchness}\label{ss:alg:branch}
Let $\AA$ be a recurrent transitive algebra.  We say that $\AA$ is
\begin{description}
  \def\makelabel#1{\textbf{\emph{#1}},}
\item[weakly branch] if for every $v\in X^*$, writing $|v|=n$, we have
  $\psi^n(\AA)\cap(\AA\otimes(v\otimes v^\circ))\neq\{0\}$, where
  $v\otimes v^\circ$ is the rank-$1$ projection on $\Bbbk v\le \Bbbk
  X^n$;
\item[weakly regular branch] there exists a non-trivial ideal
  $\KK\triangleleft \AA$, called the \emph{branching ideal}, with
  $M_X(\KK)\le\psi(\KK)$;
\item[branch] if for all $n\in\N$ the ideal $\langle
  \psi^n(\AA)\cap(\AA\otimes(v\otimes v^\circ)):v\in X^n\rangle$ has
  finite codimension in $\psi^n(\AA)$;
\item[regular branch] if there exists a finite-codimension ideal
  $\KK\triangleleft\AA$ with $M_X(\KK)\le\psi(\KK)$.
\end{description}

\begin{prop}\label{prop:alg:b}
  Let $\AA$ be an art algebra. Then it is infinite-dimensional.

  If $\AA$ is regular branch, then it is weakly regular branch and
  branch; if it is branch, then it is weakly branch. If it is weakly
  regular branch, then it is weakly branch.
\end{prop}
\begin{proof}
  Let $\AA$ be an art algebra; then it is unital. By assumption, the
  map $\psi_{x,y}:a\mapsto \psi(a)_{x,y}$ is onto.  Choose any $x\neq
  y$; then since $\psi(1)_{x,y}=0$, so $\psi_{x,y}$ is not one-to-one.
  It follows that $\AA$ is infinite-dimensional.

  Let now $\AA$ be regular branch, with branching ideal $\KK$.
  Since $\AA$ is infinite-dimensional, $\KK\neq0$, so $\AA$ is
  regular weakly branch. Assume now only that $\KK$ is non-trivial,
  and let $v\in X^*$ be any vertex. Since
  $M_{X^n}(\KK)\le\psi^n(\AA)$, we may take any $a\neq0$ in $\KK$ and
  consider the element $\psi^{-n}(a\otimes(v\otimes v^\circ))\neq0$ in
  $\AA$. This shows that $\AA$ is weakly branch. The other
  implications are of the same nature.
\end{proof}

The choice of $v$ in the definition of weakly branch algebra may have seemed
artificial; the following more general notion is equivalent:
\begin{lem}\label{lem:wbgen}
  Let $\AA$ be a weakly branch algebra. Then for any $n\in\N$ and any
  $\xi,\eta\in\Bbbk X^n$ there exists $a\neq0$ in $\AA$ with
  $(1-P_\xi)(\psi^na)=0=(\psi^na)(1-P_\eta)$, where $P_\xi,P_\eta\in M_{X^n}$
  denote respectively the projectors on $\xi,\eta$.
\end{lem}
\begin{proof}
  The weakly branch condition amounts to the lemma for $\xi=\eta$ a basis
  vector (element of $X^n$) of $\Bbbk X^n$. Write in full generality
  $\xi=\sum\xi_v v$ and $\eta=\sum\eta_v v$, the sums running over $v\in X^n$.
  Fix $w\in X^n$ and choose $b\neq0$ with $b\otimes(w\otimes
  w^\circ)\in\psi^n(\AA)$. For all $v,w\in X^n$ choose $c_{v,w}$ with
  $v\psi^n(c_{v,w})=w$; this is possible because projection on the $(v,w)$
  entry is a surjective map: $\AA\to\AA$.  Finally set
  \[a=\sum_{v,w\in X^n} \xi_v c_{v,v_0}b c_{v_0,w}\eta_w.\]
\end{proof}

\subsection{Hausdorff dimension}
Let $\AA$ be an art algebra. For every $n$, it has a representation
$\pi^n=\epsilon\psi^n:\AA\to M_{X^n}(\Bbbk)$. We define the
\emph{Hausdorff dimension} of $\AA$ as
\[\Hdim(\AA)=\liminf_{n\to\infty}\frac{\dim\pi^n(\AA)}{\dim M_{X^n}}.\]

Let us compute the Hausdorff dimension of the tree closure
$\overline\PP$ defined in~\ref{sss:treeclosure}. There,
$\pi^n(\overline\PP)$ is none other than $\PP_n$. Let
$\varpi_n=\ker\zeta_n$ denote the augmentation ideal of $\PP_n$. Then,
as a vector space, $\PP_{n+1}=M_X\otimes\varpi_n\oplus\PP$,
so
\[\dim\PP_{n+1}=\dim(M_X)(\dim\PP_n-1)+\dim\PP.\]
It follows that
\[\dim\PP_n=\frac{\dim\PP-1}{\dim M_X-1}(\dim M_X)^n + \frac{\dim M_X-\dim\PP}{\dim M_X-1},\]
and since $\dim\PP_0=1$ we have
\[\Hdim(\overline\PP)=\frac{\dim\PP-1}{\#X^2-1}.\]

If $\AA$ is a $\PP$-art algebra, we define its \emph{relative
  Hausdorff dimension} as
\[\Hdim_\PP(\AA)=\frac{\Hdim(\AA)}{\Hdim(\overline\PP)}=\Hdim(\AA)\frac{\#X^2-1}{\dim\PP-1}.\]

The following result is an analogue of Proposition~\ref{prop:gr:hd}, and
is proven along the same lines:
\begin{prop}
  Let $\AA$ be a regular branch $\PP$-art algebra. Then $\Hdim_\PP\AA$
  is a rational number in $(0,1]$.
\end{prop}
\begin{proof}
  Let $\AA$ have branching ideal $\KK$, and for all $n\in\N$ set
  $\AA_n=\pi^n(\AA)$. Let $M\in\N$ be large enough so that
  $\AA/\psi^{-2}M_{X^2}(\KK)$ maps isomorphically into
  $\AA_M/\pi^{M-2}(\KK)$. We then have, for all $n\ge M$,
  \begin{align*}
    \dim \AA_n&=\dim(\AA/\KK)+\dim\pi^n(\KK)\\
    &=\dim(\AA/\KK)+\dim(\psi \KK/M_X(\KK))+\#X^2\dim\pi^{n-1}(\KK)\\
    &=(1-\#X^2)\dim(\AA/\KK)+\dim(\psi \KK/M_X(\KK))+\#X^2\dim \AA_{n-1}.
  \end{align*}
  We write $\dim \AA_n=\alpha\#X^{2n}+\beta$, for some $\alpha,\beta$ to
  be determined; we have
  \[\alpha\#X^{2n}+\beta=(1-\#X^2)\dim(\AA/\KK)+\dim(\psi \KK/M_X(\KK))+\#X^2(\alpha\#X^{2(n-1)}+\beta),\]
  so $\beta=\dim(\AA/\KK)-\dim(\psi \KK/M_X(\KK))/(\#X^2-1)$. Then 
  set $\alpha=(\dim \AA_M-\beta)/\#X^{2M}$. We have solved the
  recurrence for $\dim \AA_n$, and $\alpha>0$ because $\AA_n$ has
  unbounded dimension, since $\AA$ is infinite-dimensional by
  Proposition~\ref{prop:alg:b}.
  
  Now it suffices to note that $\Hdim(\AA)=\alpha$ to obtain
  $\Hdim_\PP(\AA)>0$. Furthermore only linear equations with integer
  coefficients were involved, so $\Hdim(\AA)$, and $\Hdim_\PP(\AA)$,
  are rational.
\end{proof}  

\subsection{Tree enveloping algebras}\label{ss:thinalg}
Let $G$ be a recurrent group, acting on a tree $X^*$. We therefore
have a map $\Bbbk G\to\End(\Bbbk X^\omega)$, obtained by extending the
representation $G\to\aut X^\omega$ by $\Bbbk$-linearity to the group
algebra. We define the \emph{tree enveloping algebra} of $G$ as the image
$\AA$ of the group algebra $\Bbbk G$ in $\End(\Bbbk X^\omega)$.

This notion was introduced, slightly differently, by Said Sidki
in~\cite{sidki:primitive}; it has also appeared implicitly in various
places, notably~\cite{bartholdi-g:parabolic}
and~\cite{nekrashevych:bimodule}.

\begin{lem}\label{lem:codim}
  Let $\AA$ be a quotient of the group ring $\Bbbk G$, and let $H\le
  G$ be a subgroup. Let $\KK\le\AA$ be the right ideal generated by
  $\setsuch{h-1}{h\in H}$. Then
  \[\dim\AA/\KK\le[G:H].\]
\end{lem}
\begin{proof}
  It clearly suffices to prove the claim for $\AA=\Bbbk G$.  Let
  $n=[G:H]$ be the index of $H$ in $G$, and let $T$ be a right
  transversal of $H$ in $G$. Given $a\in\AA$, write $a=\sum a(g_i)
  g_i$ and each $g_i=h_it_i$ for some $h_i\in H,t_i\in T$. Then we
  have
  \[a = \sum a(g_i)h_it_i = \sum a(g_i)t_i + \sum a(g_i)(h_i-1)t_i,\]
  so $T$ generates $\AA/\KK$.
\end{proof}

\begin{thm}\label{thm:gp2a}
  Let $G$ be a recurrent transitive group, and let $\AA$ be its
  tree enveloping algebra.
  \begin{enumerate}
  \item $\AA$ is an art algebra.
  \item If $G$ is either a weakly branch group, a regular weakly
    branch group, a branch group, or a regular branch group, then
    $\AA$ enjoys the corresponding property.
  \end{enumerate}
\end{thm}
\begin{proof}
  Let $G$ be a recurrent transitive group, with decomposition
  $\phi:G\to G\wr\sym X$. Set $\FF=\Bbbk G$. We define $\psi:\FF\to
  M_X(\FF)$ by extending $\phi$ linearly: for $g\in G$, set
  \[\psi(g)=\sum_{x\in X}(g@x)\otimes(x^g\otimes x^\circ).\]
  We also let $\PP$ be the image of $\Bbbk\sym X$ in $M_X$; since
  $\sym X$ is $2$-transitive, $\PP\cong M_{\#X-1}\oplus\Bbbk$.

  By Proposition~\ref{prop:alg:unique} there is a unique image of
  $\FF$ that is an art subalgebra of $\FF$, and by construction this
  image is $\AA$.
  
  Assume that $G$ is regular branch, with branching subgroup $K$. Set
  \[\KK = \langle k-1:\,k\in K\rangle.\]
  Then $\KK$ is an ideal in $\AA$, of finite codimension by
  Lemma~\ref{lem:codim}. Since $x*K\le\phi(K)$, we have
  $\KK\otimes(x\otimes x^\circ)\le\psi(\KK)$ for all $x\in X$, and
  since $\AA$ is transitive we get $M_X(\KK)\le\psi(\KK)$, so $\AA$ is
  regular branch.
  
  Next, assume $G$ is weakly branch, and pick $v\in X^n$. There exists
  $1\neq g\in G$ with $g|_{X^\omega\setminus vX^\omega}=1$, say
  $g=v*h$. Then $g-1\neq 0$, and $0\neq\psi^n(g-1)=(h-1)(v\otimes
  v^\circ)\in \AA\otimes(v\otimes v^\circ)$, proving that $\AA$ is
  weakly branch.  The other implications are proven similarly.
\end{proof}

We note that the tree enveloping algebra corresponding to the odometer
(\S\ref{sss:odo}) or the lamplighter group (\S\ref{sss:ll})
are isomorphic to their respective group ring. Indeed these groups
have a free orbit in their action on $X^\omega$. Branch groups are at
the extreme opposite, as we will see below.

\begin{question}
  If $\AA$ is the tree enveloping algebra of a branch group $G$, does
  $\Hdim(G)>0$ imply $\Hdim(\AA)>0$? do we even have
  $\Hdim_\PP(\AA)\ge\Hdim_p(G)$ for $G\le W_p$?
\end{question}

\subsubsection{Algebraic Properties}
Recall that an algebra $\AA$ is \emph{just-infinite} if $\AA$ is
infinite-dimensional, and all proper quotients of $\AA$ are
finite-dimensional (or, equivalently, all non-trivial ideals in $\AA$
have finite codimension). The \emph{core} of a right ideal $\KK\le\AA$
is the maximal $2$-sided ideal contained in $\KK$. The \emph{Jacobson
  radical} $\jac\AA$ is the intersection of the maximal right ideals
of $\AA$. The \emph{upper nil radical} $\unil\AA$ is the sum of all
nil ideals of $\AA$.

An algebra $\AA$ is \emph{prime} if, given two non-zero ideals
$\II,\JJ\le \AA$, we have $\II\JJ\neq0$.  It is \emph{primitive} if it
has a faithful, irreducible module, or equivalently a maximal right
ideal with trivial core. It is \emph{semiprimitive}\footnote{aka
  \emph{$J$-semisimple}} if its Jacobson radical is trivial.

\begin{lem}\label{lem:kk2}
  Let $G$ be a regular branch group, with branching subgroup $K$. Let
  $\AA$ be its tree enveloping algebra, with branching ideal $\KK$. If either
  $K/[K,K]$ is finite, or $G$ is finitely generated, then $\KK/\KK^2$
  is finite-dimensional.
\end{lem}
\begin{proof}
  Consider $\KK=\langle k-1:\,k\in K\rangle\le\Bbbk G$. Then given
  $k_1,k_2\in K$ we have
  \[[k_1,k_2]-1=k_1^{-1}k_2^{-1}\big((k_1-1)(k_2-1)-(k_2-1)(k_1-1)\big)\in\KK^2,\]
  so $\KK^2$ contains $[K,K]-1$. This holds \emph{a fortiori} in
  $\AA$, so if $K/[K,K]$ is finite the result follows from
  Lemma~\ref{lem:codim}.

  If $G$ is finitely generated, then $\AA$ is also finitely generated,
  so all its finite-codimension subrings are also finitely
  generated~\cite{lewin:fgrings}.  In particular $\KK/\KK^2$ is
  finite-dimensional.
\end{proof}

\begin{thm}\label{thm:alg:ji}
  Let $\AA$ be a regular branch tree enveloping algebra. Then any ideal
  $\JJ\le\AA$ contains $M_{X^n}(\KK^2)$ for some large enough
  $n\in\N$.

  In particular, if $\KK/\KK^2$ is finite-dimensional, then $\AA$ is
  just-infinite, and if $\KK^4\neq0$, then $\AA$ is prime.
\end{thm}
\begin{proof}
  Assume $\AA$ is the tree enveloping algebra of the group $G$. Let
  $\JJ$ be a non-trivial ideal of $\AA$, and chose any non-zero
  $a\in\JJ$. Then $a=\sum a(g_i)g_i$, and the finitely many $g_i$ in
  the support of $a$ all act differently on $X^*$. The entries of
  $\psi^n(a)$, for large enough $n$, are therefore monomial; more
  precisely, there exist $v,w\in X^n$ such that the $(v,w)$ entry of
  $\psi^n(a)$, call it $b$, is in $G$, and therefore is invertible.

  Since $\JJ$ is an ideal, we have for any $v',w'\in X^n$
  \[(\KK\otimes(v'\otimes v^\circ))a(\KK\otimes(w\otimes
  (w')^\circ))=(\KK b\KK)\otimes(v'\otimes(w')^\circ)\le\JJ.\]
  It follows that $\JJ$ contains $M_{X^n}(\KK^2)$, which by assumption
  is cofinite-dimensional.

  Assume now that $\JJ,\JJ'$ are two non-zero ideals of $\AA$. By the
  above, there are $n,n'\in\N$ such that $\JJ$ contains
  $M_{X^n}(\KK^2)$ and $\JJ'$ contains $M_{X^{n'}}(\KK^2)$. For $m$
  larger than $\max\{n,n'\}$ we then have $0\neq
  M_{X^m}(\KK^4)\le\JJ\JJ'$.
\end{proof}

Recall also that an algebra $\AA$ is PI (``Polynomial Identity'') if there
exists $w\neq0$ in the free associative algebra $\Bbbk\{v_1,\dots,v_k\}$ such
that $w(a_1,\dots,a_k)=0$ for all $a_i\in\AA$. The following result is
analogous to~\ref{cor:nolaw}:
\begin{thm}\label{thm:nopi}
  Let $\AA$ be a weakly branch art algebra. Then it is not PI.
\end{thm}

We prove the theorem using the following result, which may be of independent
interest. Let $\AA$ be an algebra acting faithfully on a vector space $V$. We
say that $\AA$ \emph{separates $V$} if for every finite-dimensional subspace
$Y$ of $V$ and any $\xi\not\in Y$ there exists $a\in\AA$ with $Ya=0$ and
$\xi a\not\in\langle Y,\xi\rangle$.
\begin{prop}\label{prop:separates}
  Let $\AA$ be an algebra separating a vector space $V$. Then $\AA$ is not PI.
\end{prop}
\begin{proof}
  Let $P\in\Bbbk\{v_1,\dots,v_k\}$ be a non-commutative polynomial. We will
  find $a_1,\dots,a_k\in\AA$ and $\eta\in V$ such that $\eta
  P(a_1,\dots,a_k)\neq0$. We actually will show more, by induction: let
  $X_0\subset\{v_1,\dots,v_k\}^*$ be the set of monomials, without their
  coefficients, appearing in $P$, and let $X$ be the set of prefixes of words
  in $X_0$. For any $\eta\neq0\in V$, we construct $(a_1,\dots,a_k)\in\AA^k$
  such that $\setsuch{\eta x(a)}{x\in X}$ is an independent family.  It then of
  course follows that $\eta P(a)\neq0$.
  
  The induction starts with $X=\{1\}$. Then any $\eta\neq0$ will do. Let now
  $X$ contain at least two elements, and let $y=v_p\dots v_qv_r$ be a longest
  element of $X$. By induction, there exists $a\in\AA^k$ such that
  $Y_0=\setsuch{\eta x(a)}{x\in X\setminus\{y\}}$ is an independent family. If
  $\eta y(a)$ is linearly independent from $Y_0$, we have nothing to do.
  Otherwise, take $\xi=\eta(v_p\dots v_q)(a)$ and $Y=Y_0\setminus\{\xi\}$.
  Since $V$ is separated by $\AA$, there exists $b\in\AA$ with $Yb=0$ and $\xi
  b\not\in\langle Y,\xi\rangle$. Set $a'_i=a_i$ for $i\neq r$, and
  $a'_r=a_r+b$. Then $\setsuch{\eta x(a')}{x\in X}$ is an independent family.
\end{proof}

\begin{proof}[Proof of Theorem~\ref{thm:nopi}]
  The algebra $\AA$ is a subalgebra of $\overline\PP$, which by definition is a
  subalgebra of $\varprojlim M_X^{\otimes n}$. We may therefore assume that
  $\AA$ is a subalgebra of $\End(V)$ for the vector space $V=\varprojlim\Bbbk
  X^n$.
  
  Let $Y$ be a finite-dimensional subspace of $V$, and let $\xi\not\in V$ be
  any. Let $\pi^n$ be the projection $V\to\Bbbk X^n$. Since $Y$ is a closed
  subspace, there exists $n\in\N$ such that $v=\pi^n(\xi)\not\in\pi^n(Y)$, and
  furthermore such that there is also $w\in\pi^n(V)$ linearly independent from
  $v$ and $\pi^n(Y)$. By Lemma~\ref{lem:wbgen} there exists $a\in\AA$ which
  annihilates $Y$ while it sends $v$ to a multiple of $w$. Consider all
  possible such $a$; if they all annihilated $\xi$, then they would also
  annihilate the orbit of $\xi$ under $P_u\AA$, where $P_u\in M_{X^n}$ denotes
  projection on $u$; since they also annihilate $V(1-P_u)$, they would all
  annihilate $V$, whence $a=0$ because the representation $V$ is assumed
  faithful. This contradicts the condition that $\AA$ is weakly branch.
  
  We may therefore apply Proposition~\ref{prop:separates} to conclude that
  $\AA$ is not PI.
\end{proof}

In analogy with Proposition~\ref{prop:wb:centerless}, we have:
\begin{prop}
  Let $\AA$ be an art algebra which is weakly regular branch, with
  branching ideal $\KK$. Assume that $\KK$ is prime. Then $Z(\AA)=1$.
\end{prop}
\begin{proof}
  Take $x\in\AA$, and assume that $x$ commutes with $\KK$; we wish to
  show that $x$ is a scalar. For that, write $\psi(x)=(x_{uv})$, and
  compute $\psi[x,y\otimes(u\otimes v)]$ for all $y\in\KK$ and $u,v\in
  X$.  This matrix vanishes except possibly in its $u$th row and $v$th
  column; the $(u,v)$-entry is $x_{uu}y-yx_{vv}$, and for $v'\neq v$
  and $u'\neq u$ the $(u,v')$-entry is $yx_{vv'}$ and the
  $(u',v)$-entry is $x_{u'u}y$.

  If all those entries are to vanish, then $x_{uv}\KK=\KK x_{uv}=0$
  for all $u\neq v$, so $x_{uv}=0$ because $\KK$ is prime. Similarly
  $x_{uu}=x_{vv}$ for all $u,v$, so $\psi(x)=x_{uu}\otimes1$ for any
  $u$. Finally $[x_{uu},\KK]=0$, so the argument can be applied to
  $x_{uu}$ to show that $\psi^n(x)$ is scalar for all $n$.

  Now if $x$ were not scalar there would be $u,v\in X^n$ for some $n$
  large enough, such that $x_{uv}\neq0$ or $x_{uu}\neq x_{vv}$.
\end{proof}

\subsubsection{Compatible filtrations}
Let $\AA$ be the tree enveloping algebra of a regular branch group
$G$. We have three descending filtrations of $\AA$ by ideals, namely
powers of the branching ideal $(\KK^n)$; powers of the augmentation
ideal $(\varpi^n)$; and $(M_{X^n}(\KK))$.
\begin{prop}\label{prop:control}
  Assume that there is an $n\in\N$ such that $M_{X^n}(\KK)$ is
  contained in $\KK^2$. Then the normal subgroups of $G$
  \emph{control} the ideals of $\AA$: given any non-zero ideal
  $\JJ\le\AA$, there exists a non-trivial normal subgroup
  $H\triangleleft G$ with $H-1\subset\JJ$.
\end{prop}
\begin{proof}
  By Theorem~\ref{thm:alg:ji}, there is $n\in\N$ such that $\JJ$
  contains $M_{X^{n-1}}(\KK^2)$, so contains $M_{X^n}(\KK)$. Set
  $H=\phi^{-n}(K^{X^n})$; then $\JJ$ contains $H-1$.
\end{proof}
\begin{cor}
  Assume that there is an $n\in\N$ such that $M_{X^n}(\KK)$ is
  contained in $\KK^2$. Then $\AA$ is just-infinite and prime.
\end{cor}

Proposition~\ref{prop:control} may be used to obtain some information
on the Jacobson radical of $\AA$:
\begin{lem}[\cite{sidki:primitive}*{Corollary~4.4.3}]\label{lem:jacalternative}
  Let $\Bbbk$ be a field of characteristic $p$; let $G$ be a
  just-infinite-$p$ group (i.e. an infinite group all of whose proper
  quotients are finite $p$-groups), and let $\AA$ be a quotient of
  $\Bbbk G$. Assume that normal subgroups of $G$ control ideals of
  $\AA$.  Then either $\jac\AA=0$ or $\jac\AA=\varpi$.
\end{lem}
\begin{proof}
  $\jac\AA\le\varpi$ since $\varpi$ is a maximal right ideal. If
  $\jac\AA\neq0$, then there is a non-trivial $H\triangleleft G$ with
  $H-1\subset\jac\AA$. Since $G/H$ is a finite $p$-group, $\AA/\jac\AA$ is a
  nilpotent algebra, so is $0$, and $\jac\AA=\varpi$.
\end{proof}

This in turn gives control on representations of $\AA$, by the following
result due to Farkas and Small:
\begin{prop}[\cite{farkas-small:fd}]
  Let $\AA$ be a just-infinite, semiprimitive, finitely generated
  $\Bbbk$-algebra over an uncountable field $\Bbbk$. Then either $\AA$ is
  primitive, or $\AA$ satisfies a polynomial identity.
\end{prop}
Since weakly branch art algebras satisfy no polynomial identity
(Theorem~\ref{thm:nopi}), they admit irreducible faithful
representations as soon as they are semiprimitive.

The following are well known:
\begin{prop}[\cite{lam:ncr}*{Chapter~4}]\label{prop:lam}
  \begin{itemize}
  \item If $\AA$ is a just-infinite $\Bbbk$-algebra and contains a
    transcendental element, then $\AA$ has no non-trivial nil ideal.
  \item If $\jac\AA$ is algebraic, then it is nil.
  \item If $\AA$ is countably generated and $\Bbbk$ is uncountable,
    then the Jacobson radical $\jac\AA$ is nil.
  \item If $x\in\AA$ is transcendental and $\Bbbk$ is uncountable,
    then there exists $\alpha\in\Bbbk$ with $1-\alpha x$ not
    left-invertible.
  \end{itemize}
\end{prop}

Agata Smoktunowicz has been kind enough to explain the following
connection to me:
\begin{cor}\label{cor:ji=>p}
  If $\AA$ is just-infinite, finitely generated over an uncountable
  field $\Bbbk$, and contains a transcendental element, then $\AA$ is
  primitive.
\end{cor}

\subsubsection{The tree enveloping algebra of $\overline P$}
Consider as in \S\ref{ss:gr:tc} a subgroup $P$ of $\sym X$, and its
tree closure $\overline P\le\aut(X^*)$. It is regular branch, with
branching subgroup $\overline P$.
\begin{prop}
  Let $\AA$ be the tree enveloping algebra of $\overline P$, and let $\PP$ be
  the image in $M_X$ of $\Bbbk P$. Then $\AA=\overline\PP$.
\end{prop}
\begin{proof}
  Since $\AA\le\overline\PP$, it suffices to show that the natural map
  $\Bbbk\overline P\to\PP_n$ is onto for every $n$. Let $\varpi$
  denote the augmentation ideal of $\Bbbk\overline P$; then
  \[\psi(\Bbbk\overline P)=M_X(\varpi)+1\otimes\PP,\]
  and therefore $\psi^n(\Bbbk\overline
  P)=M_{X^n}(\varpi)+1\otimes\PP_n$, and the result follows.
\end{proof}

The algebra $\overline\PP$ can be defined in a different way,
following~\cite{sidki:primitive}. The group $\overline P$ is a
profinite (compact, totally disconnected) group, and therefore
$\Bbbk\overline P$ is a topological ring. Consider the ideal
\begin{equation}\label{eq:thinideal}
  \JJ=\langle(v*g-1)(w*h-1):\,v\neq w\in X^n\text{ for some }n;\,g,h\in\overline P\rangle
\end{equation}
in $\Bbbk\overline P$.  On the one hand, $\JJ$ has trivial image in
$\overline\PP$, since in $\psi^n(v*g-1)$ and $\psi^n(w*h-1)$ are
diagonal matrices with a single non-zero entry, in different
coordinates $v,w$. On the other hand, all relations in the matrix ring
$M_{X^n}(\Bbbk\overline P)$ can be reduced to these. It follows that
$\overline\PP$ equals $\Bbbk\overline P/\overline\JJ$, where
$\overline\JJ$, the ``thinning ideal'', denotes the
closure\footnote{note that~\cite{sidki:primitive} does not mention
  this closure, although it is essential.} of $\JJ$ in the topological
ring $\Bbbk\overline P$. More details appear in
\S\ref{ss:thinning}.

For any recurrent group $G$, we may now consider $G$ as a subgroup of
some $\overline P$, and therefore $\Bbbk G$ is a subalgebra of
$\Bbbk\overline P$. The tree enveloping ring of $\Bbbk G$ is then $\Bbbk
G/(\Bbbk G\cap\overline\JJ)$. This was the original definition of
tree enveloping rings.

\subsection{Lie algebras}
In this subsection, we let $p$ be a prime, $\Bbbk=\F$, and fix
$X=\{1,\dots,p\}$.  Let $G$ be a recurrent subgroup of $W_p$, with
decomposition $\phi:G\to G\wr C_p$ where $C_p$ is the cyclic subgroup
of $\sym X$ generated by $(1,2,\dots,p)$.  We define the
\emph{dimension series} $(G_n)$ of $G$ by $G_1=G$, and
\[G_n = \big\langle [g,h]k^p:\,g\in G,\,h\in G_{n-1},\,k\in G_{\lceil n/p\rceil}\big\rangle.\]
Since $G$ is residually-$p$, we have $\bigcap G_n=\{1\}$.

The quotient $G_n/G_{n+1}$ is an $\F$-vector space, and we form the
``graded group''
\[\gr G = \bigoplus_{n\ge1}G_n/G_{n+1}.\]
Multiplication and commutation in $G$ endows $\gr G$ with the
structure of a graded Lie algebra over $\F$, and $x\mapsto x^p$
induces a Frobenius map on $\gr G$, turning it into a \emph{restricted}
Lie algebra.

The dimension series of $G$ can be alternately described, using the
augmentation ideal $\varpi$ of $\F G$, as
\[G_n = \setsuch{g\in G}{g-1\in\varpi^n}.\]
Furthermore, consider the graded algebra $\gr\F
G=\bigoplus_{n\ge0}\varpi^n/\varpi^{n+1}$ associated to the descending
filtration $(\varpi^n)$ of $\F G$. Then
\begin{prop}[Lazard~\cite{lazard:nilp}*{Th\'eor\`eme 6.10}; Quillen~\cite{quillen:ab}]\label{thm:quillen}
  $\gr\F G$ is the restricted enveloping algebra of $\gr G$.
\end{prop}

\subsubsection{Graded tree enveloping algebras}\label{sss:gtea}
Let $\AA$ be the tree enveloping algebra of the regular branch group $G$, and
assume that $\AA$ is a graded algebra with respect to the filtration
$(\varpi^n)$. Then $\gr G$ embeds isomorphically in $\AA$.
\begin{prop}\label{prop:lieembed}
  Assume that $\AA$ is a quotient of $\gr\F G$. Then the natural map
  $\gr G\hookrightarrow\gr\F G$ induces an embedding $\gr
  G\hookrightarrow\AA$.
\end{prop}
\begin{proof}
  Let $a\in\gr G$ be such that its image in $\AA$ is trivial. Then,
  since $\AA$ is graded, all the homogeneous components of $a$ are
  trivial. But these homogeneous components belong to quotients
  $G_n/G_{n+1}$ along the dimension series of $G$, and since
  $G\hookrightarrow\AA$, they must be trivial in $G_n/G_{n+1}$. We
  deduce $a=0$.
\end{proof}

If we forget for a moment the distinction between $\Bbbk G$ and
$\gr\Bbbk G$, Proposition~\ref{prop:lieembed} can be made more
conceptual, by returning to the ``thinning process'' described
after~\eqref{eq:thinideal}: assume $G$ factors as $A\times B$. Then
$\Bbbk G=\Bbbk A\otimes\Bbbk B$, and the ``thinning'' process maps
$\Bbbk G$ to
\[\Bbbk G/\JJ=(\Bbbk A\oplus\Bbbk B)/\{(1,0)=(0,1)\},\]
with $\JJ=\varpi(\Bbbk A)\otimes\varpi(\Bbbk B)$.  We have $\gr
A\subset\Bbbk A$ and $\Bbbk B\subset\Bbbk B$ and $\gr G=\gr A\oplus\gr
B\subset\Bbbk G/\JJ$. It is in this sense that thinning ``respects''
Lie elements. More details are given in \S\ref{ss:thinning}.

Proposition~\ref{prop:lieembed} applies in particular to the group
$\overline P$ and its tree enveloping algebra $\overline\PP$. This points out
the recursive structure of $\gr\overline P$, as described
in~\cite{bartholdi:lcs}.

\subsection{Gelfand-Kirillov dimension}
Let $\AA$ be an algebra (not necessarily associative), with an
ascending filtration $(\FF_n)_{n\in\Z}$ by finite-dimensional
subspaces. Assume $\FF_n=0$ for negative $n$. Then the
\emph{Hilbert-Poincar\'e} series of $\AA$ is the formal power series
\[\Phi_\AA(t)=\sum_{n=0}^\infty a_nt^n=\sum_{n\ge0}\dim(\FF_n/\FF_{n-1})t^n.\]
In particular, if $\AA$ is generated by a finite set $S$, it has a
standard filtration defined as follows: $\FF_n$ is the linear span of
all at-most-$n$-fold products $s_1\dots s_k$ for all $k\le n$, in any
order (if $\AA$ is not associative).

If $\AA=\bigoplus_{n\ge0}\AA_n$ is graded, we naturally filter $\AA$
by setting $\FF_n=\AA_0+\dots+\AA_n$.

If $\dim\FF_n$ grows polynomially, i.e. $p_1(n)\le\dim\FF_n\le p_2(n)$
for polynomials $p_1,p_2$ of same degree, then $\AA$ has
\emph{polynomial growth}. More generally, if $\dim\FF_n$ is bounded
from above by a polynomial, the (lower) \emph{Gelfand-Kirillov}
dimension of $\AA$ is defined as
\[\GKdim(\AA)=\liminf_{n\to\infty}\frac{\log\dim\FF_n}{\log n}.\]
If $\AA$ is finitely generated and $\FF_n$ is the span of
at-most-$n$-fold products of generators, then this limit does not
depend on the choice of finite generating set.

If $\AA$ is finitely generated and either Lie or associative, then the
coefficients $a_n$ may not grow faster than exponentially. A wide
variety of intermediate types of growth patterns have been studied by
Victor Petrogradsky~\cites{petrogradsky:growth,petrogradsky:growth2}.

Let $G$ be a group, with Lie algebra $\gr G$.  Then the
Poincar\'e-Birkhoff-Witt Theorem gives a basis of $\gr\F G$ consisting
of monomials over a basis of $\gr G$, with exponents at most $p-1$. As
a consequence, we have the
\begin{prop}[Jennings~\cite{jennings:gpring}]\label{prop:sumprod}
  Let $G$ be a group with dimension series $(G_n)$, and set
  $\ell_n=\dim_{\F}(G_n/G_{n+1})$. Then
  \[\Phi_{\gr\F G}(t)=\prod_{n=1}^\infty\left(\frac{1-t^{pn}}{1-t^n}\right)^{\ell_n}.\]
\end{prop}

Approximations from analytic number theory~\cite{li:nt} and complex
analysis give then the
\begin{prop}[\cite{petrogradsky:polynilpotent}, Theorem~2.1]\label{prop:grlgr}
  With the notation above for $\ell_n$, and
  $a_n=\dim\varpi^n/\varpi^{n+1}$, we have
  \begin{enumerate}
  \item $\{a_n\}$ grows exponentially if and only if $\{\ell_n\}$ does, and
    we have
    \[\limsup_{n\to\infty}\frac{\ln \ell_n}n
    =\limsup_{n\to\infty}\frac{\ln a_n}n.\]
  \item If $\ell_n\sim n^d$, then $a_n\sim e^{n^{(d+1)/(d+2)}}$.
  \end{enumerate}
\end{prop}  

A lower bound on the growth of a group $G$ may be obtained from the
growth of $\overline{\F G}$:
\begin{prop}[\cite{grigorchuk:hp}, Lemma~8]\label{prop:gpgr}
  Let $G$ be a group generated by a finite set $S$, and let $f(n)$ be
  its growth function. Then
  \[f(n)\ge\dim(\varpi^n/\varpi^{n+1})\text{ for all }n\in\N.\]
  It follows that if $\gr G$ has Gelfand-Kirillov dimension $d$, then
  $G$ has growth at least $\exp(n^{(d+1)/(d+2)})$.
\end{prop}
It follows that a non-nilpotent residually-$p$ group has growth at
least $\exp(\sqrt n)$. It also follows that $1$-relator groups that
are not virtually abelian have exponential
growth~\cite{ceccherini-g:unifexpo}.

\begin{thm}\label{thm:gkdim}
  Let $G$ be a contracting group in the sense of
  \S\ref{sss:contraction}, acting on the tree $X^*$. Let $\AA$ be its
  tree enveloping algebra. Then $\AA$ has Gelfand-Kirillov dimension
  \begin{equation}\label{eq:gkdim}
    \GKdim(\AA)\le2n\frac{\log{\#X}}{-\log\lambda};
  \end{equation}
  in particular, if $(X^\omega,*)$ has polynomial growth of degree
  $d$, then $\AA$ has Gelfand-Kirillov dimension at most $2d$.
\end{thm}
\begin{proof}
  Let $S$ be the chosen generating set of $G$, and write
  $f(r)=\dim_\Bbbk(\Bbbk S^r)$. Then by contraction
  \[\Bbbk S^r\subset M_{X^n}(\Bbbk S^{\lambda r+K}),\]
  so $f(r)\le \#X^{2n}f(\lambda r+K)$. It follows that $\log f(r)/\log
  r$ converges to the value claimed in~\eqref{eq:gkdim}.

  The last remark follows immediately from
  Proposition~\ref{prop:contr}.
\end{proof}
\begin{question}
  Assume furthermore that $G$ is branch. Do we then have equality
  in~\eqref{eq:gkdim}?
\end{question}

\section{Examples of Tree Enveloping Algebras}\label{sec:examples}
We describe here in more detail some tree enveloping algebras. Most of the
results we obtain concern the Grigorchuk group. They are modeled on the
following result. Said Sidki considers in~\cite{sidki:primitive} the
tree enveloping algebra $\AA$ of the Gupta-Sidki group $\GS$ of
\S\ref{sss:gs}, over the field $\F[3]$. He shows:

\begin{thm}\label{thm:sidki}
  \begin{enumerate}
  \item The group $\GS$ and the polynomial ring $\F[3][t]$ embed in
    $\AA$;
  \item The algebra $\AA$ is just-infinite, prime, and primitive.
  \end{enumerate}
\end{thm}

\subsection{The ``thinning process''}\label{ss:thinning}
We recall and generalize the original construction of $\AA$, since it
is relevant to \S\ref{sss:gtea}. Let $G\hookrightarrow G\wr
P$ be a recurrent group, with $P\le\sym X$. Let $\FF=\Bbbk G$ be its
group algebra. Then we have a natural map
\[\FF\hookrightarrow \FF^{\otimes X}\rtimes P=\FF^{\otimes
  X}\underline\otimes\Bbbk P,
\]
where $\AA\rtimes P$ designates the crossed product algebra; the
$\underline\otimes$ indicates the tensor product as vector spaces, with
multiplication
\[(1^{\otimes X}\underline\otimes\pi)(g_1\otimes\dots\otimes
g_q\underline\otimes1)=(g_{1^\pi}\otimes\dots\otimes
g_{q^\pi}\underline\otimes1)\otimes(1^{\otimes X}\underline\otimes\pi)
\]
for all $g_1,\dots,g_q\in G$ and $\pi\in P$.

We wish to construct a quotient of $\FF$ which still contains a copy
of $G$. For this, let $\varpi_i$ denote, for all $i\in X$, the
augmentation ideal of the subalgebra
$\Bbbk\otimes\dots\otimes\FF\otimes\dots\otimes\Bbbk\cong\FF$, with
the `$\FF$' in position $i$; and let $\II_i$ denote the ideal in
$\Bbbk P$ generated by $\setsuch{\pi-1}{i^\pi=i}$. Set then
\[\JJ = \sum_{i\neq j\in X}\varpi_i\varpi_j\otimes\Bbbk P + \sum_{i\in
  X}\varpi\otimes\II_i + \bigcap_{i\in X}\Bbbk\otimes\II_i.\]
\begin{lem}[\cite{sidki:primitive}*{\S 3.2}]
  $\FF/\JJ\cong M_X(\FF)$.
\end{lem}
This process can then be iterated, by thinning the `$\FF$' on the
right-hand side of the above; the limit coincides with the tree
enveloping ring of $G$.

\subsection{The Grigorchuk group}
From now on, we restrict to the Grigorchuk group $\Gg$ defined in
\S\ref{sss:gr:gg}. There are two main cases to consider,
depending on the characteristic of $\Bbbk$: tame ($\neq2$) or wild
($=2$).

We begin by some general considerations. As generating set of $\Gg$ we
always choose $S=\{a,b,c,d\}$, and we may again choose $S$ as
generating set of its tree enveloping algebra $\AA$.

Since $\Gg$'s decomposition is $\phi:\Gg\mapsto\Gg\wr C_2$, the ring
$\PP$ is the linear envelope of the representation of $C_2$ on two
points, i.e.\ the group ring of $C_2$:
\[\PP=\left.\left\{\mat\alpha\beta\beta\alpha\right|\,\alpha,\beta\in\Bbbk\right\}\cong\Bbbk[\Z/2].\]
If $\Bbbk$ has characteristic $2$, this is the nilpotent ring
$\Bbbk[t]/(t^2)$; in tame characteristic, $\PP=\Bbbk\oplus\Bbbk$.

Following Theorem~\ref{thm:gp2a}, we may rewrite $\Gg$'s
decomposition~\eqref{eq:grgr:decomp} as a map $\psi:\AA\to M_2(\AA)$:
\begin{equation}\label{eq:gr:decomp}
  a\mapsto\mat 0110,\quad b\mapsto\mat a00c,\quad c\mapsto\mat
  a00d,\quad d\mapsto\mat 100b.
\end{equation}

\begin{thm}\label{thm:galg:jip}
  The algebra $\AA$ is regular branch, just-infinite, and prime.
\end{thm}
\begin{proof}
  $\AA$ is regular branch by Theorem~\ref{thm:gp2a}. By
  Lemma~\ref{lem:kk2} and Theorem~\ref{thm:alg:ji} it is
  just-infinite and prime.
\end{proof}
Ana Cristina Vieira proved in~\cite{vieira:modular}*{Corollary~4} that
$\AA$ is just-infinite if $\Bbbk=\F[2]$. Actually her arguments extend
to arbitrary characteristic, and also show that $\AA$ is prime.

\subsubsection{Characteristic $\neq2$}
In this subsection, let $\Bbbk$ be a field of characteristic $\neq2$.
\begin{prop}
  The algebra $\AA$ is semiprimitive. If furthermore $\Bbbk$ is
  uncountable, then $\AA$ is primitive.
\end{prop}
\begin{proof}
  The ring $\AA$ admits finite-dimensional quotients
  $\AA_n=\AA/\PP_n=\pi^n(\AA)$. Since $\Bbbk$ was assumed of
  characteristic $\neq2$ and $\AA_n$ is a quotient of the group
  algebra of a $2$-group, it is semisimple and therefore
  $\jac\AA\le\PP_n$ for all $n$, so $\jac\AA=0$.

  If $\Bbbk$ is uncountable, then $\AA$ is primitive
  by~\cites{amitsur:semisimple,passman:nil}.
\end{proof}
\begin{question}
  Is $\AA$ primitive for $\Bbbk=\Q$ or $\F$ with $p\neq2$?
\end{question}

\begin{prop}
  The algebra $\AA$ has relative Hausdorff dimension $\Hdim_\PP(\AA)=1$.
\end{prop}
\begin{proof}
  This is a reformulation
  of~\cite{bartholdi-g:parabolic}*{Theorem~9.7}, where the structure
  of the finite quotient $\pi^n(\AA)$ is determined for $\Bbbk=\C$.
  The result obtained was
  \[\pi^n(\AA)=\C+\bigoplus_{i=0}^{n-1} M_{2^i}(\C).\]
  It follows that $\pi^n(\AA)$ has dimension $(4^n+2)/3$. The proof
  carries to arbitrary $\Bbbk$ of characteristic $\neq2$.
\end{proof}

The algebra $\AA$ does not seem to have any natural grading; indeed if
$\varpi$ denote the augmentation ideal of $\AA$, then
$\varpi^2=\varpi$, because $\varpi$ is generated by idempotents
$\frac12(1-a),\frac12(1-b),\frac12(1-c),\frac12(1-d)$. As a side note,
the Lie powers $\varpi^{[n]}$ of $\varpi$, defined by
$\varpi^{[1]}=\varpi$ and
\[\varpi^{[n+1]}=\AA\setsuch{xy-yx}{x\in\varpi^{[n]},y\in\varpi}\AA,\]
also seem to stabilize.

The following presentation is built upon
Proposition~\ref{prop:lysenok}.  Since the proof is similar to that of
Theorem~\ref{thm:=2:pres}, we only sketch the proof.
\begin{thm}\label{thm:<>2:pres}
  Consider the endomorphism $\sigma$ of $\Bbbk\{a,b,c,d\}$ defined on
  its basis by
  \begin{equation}\label{eq:<>2sigma}
    a\mapsto aca,\quad b\mapsto d,\quad d\mapsto c,\quad c\mapsto b
  \end{equation}
  and extended by linearity. Then
  \begin{multline}
    \AA=\big\langle a,b,c,d\big|a^2=b^2=c^2=d^2=bcd=1,\\
    \sigma^n\big((d-1)a(d-1)\big)=\sigma^n\big((d-1)a(d^{acac}-1)\big)=0\;\forall n\ge0\big\rangle.\label{eq:<>2lysenok}
  \end{multline}
\end{thm}
\begin{proof}
  Let $\FF$ be the free associative algebra on $S$; define
  $\psi:\FF\to M_2(\FF)$ using formul\ae~\eqref{eq:gr:decomp}. Set
  $\JJ_0=\langle a^2-1,b^2-1,c^2-1,d^2-1,bcd-1\rangle$,
  $\JJ_{n+1}=\psi^{-1}(M_2(\JJ_n))$, and $\JJ=\bigcup_{n\ge0}\JJ_n$.
  We therefore have an algebra $\AA'=\FF/\JJ$, and since an easy check
  shows that the relations above hold in $\AA$, we have a natural map
  $\pi:\AA'\to\AA$ which is onto. We show that it is also one-to-one.

  Take $x\in\ker\pi$. Then it is a finite linear combination of words
  in $S^*$, so there exists $n\in\N$ such that all entries in
  $\psi^n(x)$ are linear combinations of words of syllable length at
  most $1$, where $a$'s and $\{b,c,d\}$'s are grouped in syllables.
  Since they must also act trivially on $\Bbbk X^\omega$, they belong
  to $\JJ_0$; so $x\in\JJ_n$.

  It remains to compute $\JJ_n$. First, $\JJ_1/\JJ_0$ is generated by
  all $(d^u-1)a(d^v-1)$ for $u,v\in\{a,b,c,d\}^*$ with an even number
  of $a$'s. It is sufficient to consider only $u=1$; and to assume
  that $v$ contains only $a$'s and $c$'s; indeed $d$'s can be
  pulled out to give a shorter relator of the form $(d-1)a(d^w-1)$,
  and $b$'s can be replaced by $c$'s by the same argument. Using the
  previous relators, we may then suppose that $v$ is of the form
  $(ac)^{2k}$.

  Next, the relators $r_k=(d-1)a(d^{(ac)^{2k}}-1)\in\JJ_1$ lift to
  generators $\sigma^n(r_k)$ of $\JJ_{n+1}/\JJ_n$.

  Finally, using the relator $\sigma(r_0)=cacac-aca$, we see that it
  is sufficient to consider the relators $\sigma^n(r_0)$ and
  $\sigma^n(r_1)$.
\end{proof}

Although we may not grade $\AA$, we may still filter it by powers of
the generating set $S$.  We give the following result with minimal
proof; it follows from arguments similar, but harder, than those in
Proposition~\ref{prop:grig:powers}.
\begin{thm}\label{thm:<>2:gkdim}
  The algebra $\AA$ has quadratic growth; therefore its
  Gelfand-Kirillov dimension is $2$.

  More precisely, set $\FF_n=\sum_{i=0}^n\Bbbk S^i$ and
  $a_n=\dim\FF_n/\FF_{n-1}$. Then
  $a_1=4,a_2=6,a_3=8,a_4=10,a_5=13,a_6=16$, and for $n\ge7$
  \begin{equation}\label{eq:dim<>2}
    a_n=\begin{cases}
      4n-\frac322^k & \text{ if }2^k\le n\le\frac542^k,\\
      3n-\frac142^k & \text{ if }\frac542^k\le n\le\frac322^k,\\
      n+\frac{11}42^k & \text{ if }\frac322^k\le n\le\frac742^k,\\
      2n+2^k & \text{ if }\frac742^k\le n\le2^{k+1}.\end{cases}
  \end{equation}
  It follows for example that, if $n$ is a power of two greater than
  $4$, then
  \[\dim\FF_n=\frac43n^2+\frac54n+\frac23.\]
\end{thm}
Note that $\AA$ has Gelfand-Kirillov dimension at most $2$, by
Theorem~\ref{thm:gkdim}; furthermore, it cannot have dimension $1$
since $\AA$ satisfies no polynomial identity by
Theorem~\ref{thm:nopi}, so by Bergman's gap
theorem~\cite{krause-l:gkdim} it has dimension $2$.

\begin{lem}\label{lem:<>2:J}
  Set $x=ab-ba$ and let $\KK=\AA x\AA$ be the branching ideal of
  $\AA$. Then $\AA/\KK$ is $6$-dimensional, and $\KK/M_2(\KK)$ is
  $20$-dimensional.
\end{lem}
\begin{proof}
  The codimension of $\KK$ is at most $16$, which is the index of $K$
  in $\Gg$. We then check $y=(1+b)(1-d)\in\KK$, because
  $y=\frac12(c-1)xay$, and we use $(d-1)a(d-1)=0$ to see that the
  codimension of $\KK$ is at most $6$, with transversal
  $\{1,a,d,ad,da,ada\}$. These elements are easily seen to be
  independent modulo $\KK$.

  The assertion on $\KK/M_2(\KK)$ has a similar proof.
\end{proof}
\begin{proof}[Sketch of the proof of Theorem~\ref{thm:<>2:gkdim}]
  The first few values of $a_n$ are computed directly.  We consider
  the filtrations $\EE_n=\FF_n\cap\KK$ and $\DD_n=\FF_n\cap M_2(\KK)$
  of $\KK$ and $M_2(\KK)$ respectively. For $n\ge3$, we have
  $\dim\FF_n/\EE_n=6$, and for $n\ge6$ we have
  $\dim\EE_n/\DD_n=20$. It follows that $a_n=\dim\DD_n/\DD_{n-1}$ for
  $n$ large enough, and we place ourselves in that situation.

  A word $w\in S^*$ is \emph{reduced} if it alternates between
  $a$-letters and $\{b,c,d\}$-letters. Every group element in $\Gg$
  can be represented by a reduced word. We construct the following
  refinement of the filtration $(\FF_n)$: we denote by $\FF_n^{at}$
  the linear span of those words $w\in S^*$ which either have length
  $\le n-1$ or are reduced, of length $n$, start in $a$, and end in
  $\{b,c,d\}$. We define similarly $\FF_n^{ta},\FF_n^{aa},\FF_n^{tt}$.
  We set $\EE_n^{at}=\FF_n^{at}\cap\KK$, and define similarly
  $\EE_n^{ta},\EE_n^{aa},\EE_n^{tt},\DD_n^{at},\DD_n^{ta},\DD_n^{aa},\DD_n^{tt}$.
  If $n$ is even, then $\DD_n=\DD_n^{at}+\DD_n^{ta}$, while if $n$ is
  odd, then $\DD_n=\DD_n^{aa}+\DD_n^{tt}$.

  The following equalities are not hard to check; the ``$\subseteq$''
  part comes from the contraction of $\Gg$'s decomposition map, and
  the ``$\supseteq$'' part comes from a construction using the
  endomorphism $\sigma$ of~\eqref{eq:<>2sigma}. For $n\ge3$. we have
  \begin{align*}
    \DD_{4n} &= \DD_{4n-1} + \DD_{4n}^{at}+\DD_{4n}^{ta}
    = \DD_{4n-1} + \mat{\EE_{2n}^{ta}}00{\EE_{2n}^{at}}+\mat{\EE_{2n}^{at}}00{\EE_{2n}^{ta}}\\
    &= \DD_{4n-1} + \mat{\EE_{2n}}00{\EE_{2n}},\\
    \intertext{and similarly}
    \DD_{4n+1} &= \DD_{4n} + \DD_{4n+1}^{aa}+\DD_{4n+1}^{tt}
    = \DD_{4n} + \mat0{\EE_{2n}^{ta}}{\EE_{2n}^{at}}0 + \mat{\EE_{2n+1}^{aa}}00{\EE_{2n+1}^{tt}},\\
    \DD_{4n+2} &= \DD_{4n+1} + \DD_{4n+2}^{at}+\DD_{4n+2}^{ta}
    = \DD_{4n+1} + \mat0{\EE_{2n+1}^{aa}}{\EE_{2n+1}^{tt}}0+\mat0{\EE_{2n+1}^{tt}}{\EE_{2n+1}^{aa}}0\\
    &= \DD_{4n+1} + \mat0{\EE_{2n+1}}{\EE_{2n+1}}0,\\
    \DD_{4n+3} &= \DD_{4n+3}^{aa}+\DD_{4n+3}^{tt}
    = \DD_{4n+2} + \mat{\EE_{2n+1}^{tt}}00{\EE_{2n+1}^{aa}} + \mat0{\EE_{2n+2}^{at}}{\EE_{2n+2}^{ta}}0.
  \end{align*}
  These equalities give
  \begin{align*}
    a_{4n}&=\dim(\DD_{4n}/\DD_{4n-1})=2\dim(\EE_{2n}/\EE_{2n-1})=2a_{2n},\\
    a_{4n+1}&=\dim(\EE_{2n}^{ta}/\EE_{2n-1})+\dim(\EE_{2n}^{at}/\EE_{2n-1})+\dim(\EE_{2n+1}^{aa}/\EE_{2n})\\
    &\kern+5cm{}+\dim(\EE_{2n+1}^{tt}/\EE_{2n}) = a_{2n}+a_{2n+1},\\
    a_{4n+2}&=2\dim(\EE_{2n+1}/\EE_{2n})=2a_{2n+1},\\
    a_{4n+3}&=\dim(\EE_{2n+1}^{tt}/\EE_{2n})+\dim(\EE_{2n+1}^{aa}/\EE_{2n})+\dim(\EE_{2n+2}^{at}/\EE_{2n+1})\\
    &\kern+5cm{}+\dim(\EE_{2n+2}^{ta}/\EE_{2n+1}) = a_{2n+1}+a_{2n+2},
  \end{align*}
  from which~\eqref{eq:dim<>2} follows.
\end{proof}

\subsection{The Grigorchuk group in characteristic $2$} If we let
$\Bbbk$ be a field of characteristic $2$, then sharper results appear.
To state them, it is better to choose another generating set for
$\AA$, and throughout this subsection we assume $S=\{A,B,C,D\}$, with
$A=a-1,B=b-1,C=c-1,D=d-1$. In that notation, the augmentation ideal
$\varpi$ of $\AA$ is generated by $S$, and $\AA$ is generated by $S$
as an algebra with one.

We first recall, in a more concrete form, the results stated above for
general $\Bbbk$.
\begin{prop}
  The algebra $\AA$ is recurrent; its decomposition map $\psi:\AA\to
  M_2(\AA)$ is given by
  \begin{equation}\label{eq:=2:decomp}
    A\mapsto\mat 1111,\quad B\mapsto\mat A00C,\quad C\mapsto\mat
    A00D,\quad D\mapsto\mat 000B.
  \end{equation}
\end{prop}
\begin{proof}
  The expression of $\psi$ follows from the definition. Upon
  inspection, one sees $1$, $B$, $C$ and $D$ in the $(2,2)$ corner as
  $\psi(A)$, $\psi(D)$, $\psi(B)$ and $\psi(C)$; then $\psi(ACA+C)$
  gives an $A$ in the $(2,2)$ corner, so projection on the $(2,2)$
  corner is onto. For the other corners, it suffices to multiply the
  above expressions by $1+A$ on the left, on the right, or on both
  sides to obtain all generators in the image of the $(i,j)$
  projection.
\end{proof}

\begin{thm}
  The relative Hausdorff dimension of $\AA$ is $\Hdim_\PP(\AA)=7/8$.
\end{thm}
\begin{proof}
  Let $\AA_n$ be the finite quotient $\pi^n(\AA)$ of $\AA$, and set
  $b_n=\dim\AA_n$. Then $b_2=8$ by direct examination, and one solves
  the recurrence, for $n\ge3$,
  \begin{multline*}
    b_{n+1}=\dim\AA_n=\dim\AA/\KK+\dim\pi^{n+1}(\KK)\\
    =6+\dim(\KK/M_2(\KK))+2^2\dim\pi^n(\KK)=6+8+4(b_n-6)
  \end{multline*}
  to $b_n=(14\cdot4^{n-2}+10)/3$. This gives $\Hdim(\AA)=14/24$, and
  $\Hdim_\PP(\AA)=7/8$.
\end{proof}

Let $H$ be the stabilizer in $\Gg$ of the infinite ray $1^\omega\in
X^\omega$; then by~\cite{bartholdi-g:parabolic} it is a weakly maximal
subgroup, i.e.\ if $H\lneqq I\le\Gg$ then $I$ has finite index in
$\Gg$.  It follows that the right ideal $\JJ=(H-1)\AA$ is a ``weakly
maximal'' right ideal, i.e.\ if $\JJ\lneqq\II\le\AA$ then $\II$ has
finite codimension in $\AA$. Since the core of $\JJ$ is trivial, it
follows that $\AA$ admits a faithful module $\AA/\JJ$ all of whose
quotients are finite.  This is none other than the original
representation on $\Bbbk X^*$.
\begin{prop}
  The ideal $\JJ$ has Gelfand-Kirillov dimension $1$; i.e. the
  dimensions of the quotients
  $\JJ\cap\varpi^n\big/\JJ\cap\varpi^{n+1}$ are bounded.
\end{prop}
\begin{proof}
  This is a reformulation of~\cite{bartholdi-g:lie}*{Lemma~5.2}, where
  the uniseriality of the modules naturally associated with $X^m$ is
  proven.
\end{proof}

From now on, we identify $\AA$ with its image in $M_2(\AA)$. We also
commit the usual crime of identifying words over $S$ with their
corresponding elements in $\AA$. Set
\begin{equation}
  \mathcal R_0=\{A^2,B^2,C^2,D^2,B+C+D,BC,CB,BD,DB,CD,DC,DAD\}.\label{eq:r0}
\end{equation}
We also set $T=\{B,C,D\}$.
\begin{lem}\label{lem:alg:rel}
  All words in $\mathcal R_0$ are trivial in $\AA$. Furthermore, the
  last relator is part of a more general pattern: $DwD$ is trivial for
  any word $w\in S^*$ with $|w|\equiv1\mod 4$.
\end{lem}
\begin{proof}
  Clearly $A^2=0$. Then $B+C+D=\smallmat000{B+C+D}$ so $B+C+D$ acts
  trivially on $\Bbbk X^\omega$ and is therefore trivial. Given any
  $x,y\in T$ we have $xy=\smallmat000{x'y'}$ for some $x',y'\in T$ and
  these are therefore also relations. Finally, let $w\in S^*$ be a
  word of length $4n+1$. Clearly, by the above, $DwD=0$ unless
  possibly if $w$ is of the form $Ax_1\dots Ax_{2n}A$ for some $x_i\in
  T$. Then $w=\smallmat{w_{11}}{w_{12}}{w_{21}}{w_{21}}$ where each
  $w_{ij}$ is a linear combination of words that either start or end
  in $T$; multiplying on both sides with $D=\smallmat000B$ therefore
  annihilates $DwD$.
\end{proof}

\subsubsection{A recursive presentation for $\AA$}
Consider the substitution $\sigma:S^*\to S^*$, defined as follows:
\[A\mapsto ACA,\quad B\mapsto D,\quad C\mapsto B,\quad D\mapsto C.\]
We say that a word $w\in S^*$ is an \emph{$A\div T$ word} if its first
letter is $A$ and its last letter is in $T$; we define similarly
$A\div A$, $T\div A$, and $T\div T$ words. A \emph{$\div A$ word} is a
word ending in $A$, and $\div T$, $A\div$ and $T\div$ words are
defined similarly.
\begin{lem}\label{lem:gr:sigma}
  Let $w\in S^*$ represent an element of $\KK$. Then in $\AA$ we have
  \begin{itemize}
  \item if $w$ is a $A\div A$ word, then $\sigma(w)=\mat wwww$;
  \item if $w$ is a $A\div T$ word, then $\sigma(w)=\mat 0w0w$;
  \item if $w$ is a $T\div A$ word, then $\sigma(w)=\mat 00ww$;
  \item if $w$ is a $T\div T$ word, then $\sigma(w)=\mat 000w$, unless
    if $w$ belongs to $\{CAC,CAD,DAC,DAD\}$, in which case
    $\sigma(w)=\mat{ADA}00w$.
  \end{itemize}
\end{lem}
Note in particular that because of the four exceptional cases for
$T\div T$ words, the map $\sigma$ does not induce an endomorphism of
$\AA$.  It seems that there does not exist a graded endomorphism
$\tau$ of $\AA$ with $\psi(\tau(w))_{2,2}=w$ for all long enough $w\in
S^*$.

\begin{proof}
  The induction starts with the words $B,CAC,CAD,DAC,DAD$. If for
  example $w$ is a $A\div T$ word, we have $\sigma(w)=\smallmat 0w0w$,
  and therefore
  \[\sigma(wA)=\mat 0w0wACA=\mat
  0w0w\mat{A+D}{A+D}{A+D}{A+D}=\mat{wA}{wA}{wA}{wA},
  \]
  where $wD=0$ because $w$ ends in a letter in $T$.
\end{proof}

\begin{prop}
  The algebra $\AA$ is regular branch.
\end{prop}
\begin{proof}
  This follows from Theorem~\ref{thm:gp2a}. Alternatively, consider
  the ideal
  \[\KK=\langle ADA,AB,BA\rangle.\]
  Compute $\dim(\AA/\KK)=6$, with $\AA=\KK\oplus\langle
  1,A,B,D,AD,DA\rangle$. Next check
  \begin{align*}
    \mat{ADA}000 &= CACAC = C(ADA)C + CA(BA)C\in \KK\\
    \mat{AB}000 &= CADA = C(ADA)\in \KK\\
    \mat{BA}000 &= ADAC = (ADA)C\in \KK,
  \end{align*}
  giving $M_2(\KK)\le \KK$. We have $\dim\KK/M_2(\KK)=8$, because
  \[\KK=M_2(\KK)\oplus\langle ADA,AB,BA,ABA,BAB,ABAB,BABA,ABABA\rangle.\]

  \noindent We may also easily check that $\KK/\KK^2$ is
  $12$-dimensional, by
  \begin{multline*}
    \KK^2=\KK\oplus\langle AB,BA,ABA,ADA,BAB,BAD,DAB,\\
    ABAD,ADAB,BADA,DABA,DABAD\rangle.
  \end{multline*}
\end{proof}

\begin{thm}\label{thm:=2:pres}
  Let $\mathcal R_0$ be as in~\eqref{eq:r0}. Then the algebra $\AA$
  admits the presentation
  \[\AA=\langle
  A,B,C,D|\,\mathcal R_0,\sigma^n(CACACAC),\sigma^n(DACACAD)\text{
    for all }n\ge0\rangle.\]
\end{thm}
\begin{cor}\label{cor:=2:graded}
  $\AA$ is graded along powers of its augmentation ideal $\varpi$.
  This grading coincides with that defined by the generating set $S$.
\end{cor}
\begin{proof}
  All relations of $\AA$ are homogeneous --- they are even all
  monomial, except for $B+C+D$.
\end{proof}
\begin{proof}[Proof of Theorem~\ref{thm:=2:pres}]
  Let $\FF$ be the free associative algebra on $S$; define
  $\psi:\FF\to M_2(\FF)$ using formul\ae~\eqref{eq:=2:decomp}. Set
  $\JJ_0=\langle\mathcal R_0\rangle$,
  $\JJ_{n+1}=\psi^{-1}(M_2(\JJ_n))$, and $\JJ=\bigcup_{n\ge0}\JJ_n$.
  We therefore have an algebra $\AA'=\FF/\JJ$, with a natural map
  $\pi:\AA'\to\AA$ which is onto. We show that it is also one-to-one.

  Take $x\in\ker\pi$. Then it is a finite linear combination of words
  in $S^*$, so there exists $n\in\N$ such that all entries in
  $\psi^n(a)$ are words in $A^*$ or $T^*$. Since they must also act
  trivially on $\Bbbk X^\omega$, they belong to $\JJ_0$; so
  $x\in\JJ_n$.

  It remains to compute $\JJ_n$. First, $\JJ_1/\JJ_0$ is generated by
  all $DwD$ with $|w|\equiv1\mod 4$, which map to $0\in\FF/\JJ_0$, and
  $CACACAC$, which maps to $DAD=0\in\FF/\JJ_0$. Using the relation
  $r_0=DAD$, we see that all $DwD$ are consequences of $r_1=DACACAD$
  and $r_2=CACACAC$. For example, $DACABAD=r_1+DACAr_0$,
  $r'_1=DABABAD=DACABAD+r_0ABAD$, and for $n\ge2$, by induction
  \begin{multline*}
    r'_n=D(AB)^{2n}AD=r'_{n-1}ABAD+r'_{n-2}ABACABAD\\
    +D(AB)^{2n-4}A\big(CABACAr_0+CAr_1+r_2AD\big)
  \end{multline*}
  
  Finally, the relations $r_1,r_2\in\JJ_1$ lift to generators
  $\sigma^n(r_1),\sigma^n(r_2)$ of $\JJ_{n+1}/\JJ_n$.
\end{proof}

\begin{prop}\label{prop:grig:powers}
  Successive powers of the augmentation ideal of $\AA$ satisfy, for
  $n\ge3$,
  \[\dim(\varpi^n/\varpi^{n+1})=\begin{cases}
    2n-\frac122^k & \text{ if }2^k\le n\le\frac322^k\\
    n+2^k & \text{ if }\frac322^k\le n\le2^{k+1}.\end{cases}
  \]
\end{prop}

It follows that, although $\Bbbk\Gg$ has large growth, namely
$\dim(\varpi^n/\varpi^{n+1})\sim\exp(\sqrt n)$ in $\Bbbk\Gg$ by
Proposition~\ref{prop:gpgr}, the growth of its quotient $\AA$ is
polynomial of degree $2$:
\begin{cor}\label{cor:=2:gkdim}
  The algebra $\AA$ has quadratic growth; therefore its
  Gelfand-Kirillov dimension is $2$, both as a graded algebra (along
  powers of $\varpi$), and as a finitely generated filtered algebra.
\end{cor}
\begin{proof}[Proof of Proposition~\ref{prop:grig:powers}]
  Assume $n\ge 3$. Then we have
  \begin{align}
    \varpi^{2n}&=\left\langle\varpi^n\otimes\mat0101,\varpi^n\otimes\mat0011\right\rangle,\label{eq:powers:1}\\
    \varpi^{2n+1}&=\left\langle\varpi^n\otimes\mat1111,\varpi^{n+1}\otimes\mat0001\right\rangle.\label{eq:powers:2}
  \end{align}
  Indeed consider a generator $w\in S^*$ of $\varpi^{2n}$. Then $w$ is
  a word of length $2n$, so is either a $A\div T$ word or a $T\div A$
  word. It follows that $\psi(w)=\smallmat 0u0u$ or $\smallmat 00uu$
  for some $u\in S^n$, and the `$\subseteq$' inclusion is shown.

  Conversely, take $u\in S^n$; if the length of $u$ is even, then $u$
  is either a $T\div A$ word or a $A\div T$ word, and set
  $w=\sigma(u)$. If $|u|$ is odd, then $u$ is either a $T\div T$ word,
  and consider $w=\sigma(u)A$ and $A\sigma(u)$, or it is a $A\div A$
  word, and set $w'=\sigma(u)$ and $w=w'$ with its first or last
  letter removed. In all cases, $w$ is a word of length $2n$, and
  $\psi(w)=\smallmat 0u0u$ or $\smallmat 00uu$, which shows the
  `$\supseteq$' inclusion. A similar argument applies
  to~\eqref{eq:powers:2}.

  Set $a_n=\dim(\varpi^n/\varpi^{n+1})$. Then it is easy to compute
  \begin{xalignat*}{2}
    \AA/\varpi&=\langle1\rangle & \text{giving }&a_0=1\\
    \varpi/\varpi^2&=\langle A,B,D\rangle & \text{giving }&a_1=3\\
    \varpi^2/\varpi^3&=\langle AB,BA,AD,DA\rangle & \text{giving }&a_2=4\\
    \varpi^3/\varpi^4&=\langle ABA,ADA,BAB,BAD,DAB\rangle & \text{giving }&a_3=5\\
    \varpi^4/\varpi^5&=\langle ABAB,ABAD,ADAB,BABA,BADA,DABA\rangle & \text{giving }&a_4=6\\
    \varpi^5/\varpi^6&=\langle ABABA,ABADA,ADABA,BABAB,\\
    &\hspace{2cm} BABAD,BADAB,DABAB,DABAD\rangle & \text{giving }&a_5=8,
  \end{xalignat*}
  and formul\ae~(\ref{eq:powers:1},\ref{eq:powers:2}) give
  \[a_{2n}=2a_n,\qquad a_{2n+1} = a_n+a_{n+1},\]
  from which the claim follows.
\end{proof}

\noindent We now show that the filtrations of $\AA$ by $(\omega^n)$,
$(\KK^n)$ and $(M_{X^n}(\KK))$ are equivalent:
\begin{prop}\label{prop:samefilt}
  For all $n\in\N$ we have
  \begin{align*}
    \varpi^{3n}\le \KK^n &\le\varpi^{2n},\\
    \varpi^{3\cdot 2^n}\le M_{X^n}(\KK) &\le\varpi^{2\cdot 2^n}.
  \end{align*}
\end{prop}
\begin{proof}
  To check the first assertion, it suffices to note that all
  non-trivial words of length $3$ in $S$, namely
  $(AB)A,ADA,(BA)B,(BA)D,D(AB)$, belong to $\KK$, while all generators
  of $\KK$ lie in $\varpi^3$.

  To check the third inclusion, take $w\in S^{3\cdot 2^n}$; then
  $\psi^n(w)\in M_{X^n}(\varpi^3)$. To check the fourth inclusion,
  take a generator $w$ of $\KK$, and consider $v=\sigma^n(w)$. Since
  $|w|\ge2$, we have $|v|\ge2\cdot2^n$ so $v\in\varpi^{2\cdot2^n}$.
\end{proof}

\subsubsection{Laurent polynomials in $\AA$}
It may seem, since $\AA$ has Gelfand-Kirillov dimension $2$, that $\Gg$
contains ``most'' of the units of $\AA$. However, $\Gg$ has infinite
index in $\AA^\times$, and contains an element of infinite order:
\begin{thm}\label{thm:=2:laurent}
  $\AA$ contains the Laurent polynomials $\Bbbk[X,X^{-1}]$.
\end{thm}
\begin{proof}
  Consider the element $X=1+A+B+AD$. It is invertible, with
  \[X^{-1}=(1+B)(1+AC)(1+ACAC)(1+A).\] Now to show that $X$ is
  transcendental, it suffices to show that $X$ has infinite order;
  indeed if $X$ were algebraic, it would generate a finite extension
  of a finite field, and therefore a finite ring; so $X$ would have
  finite order.

  Among words $w\in\{A,B,AD\}^*$, consider the set $\mathcal W$ of
  those of the form
  \[w=(AB)^{i_1}AD(AB)^{i_2}AD\dots(AB)^{i_\ell}.\]
  These are precisely the words starting by an $A$, and ending by a
  $B$ or a $D$. Define their \emph{length} and \emph{weight} as
  \[|w|=\sum_{j=1}^\ell(2i_j+2),\qquad \|w\|=\sum_{j=1}^\ell(2i_j+1).\]

  Consider the words $w_n$ defined iteratively as follows: $w_1=ADAB$,
  and $w_n=\tau(w_{n-1})$ where $\tau$ is the substitution
  $\tau(AB)=(ADAB)^3(AB)^2$, $\tau(AD)=(ADAB)^4$. Then
  \[\psi^3(w_n)=w_{n-1}\otimes\mat0101\otimes\mat1100\otimes\mat1010.\]

  Define $\sigma(n)=\frac{22\cdot 8^n-1}7$. Then $|w_n|=4\cdot8^n$ and
  $\|w_n\|=\sigma(n)$; and $w_n$ is the unique summand of
  $X^{\sigma(n)}$ in $\mathcal W$ that belongs to
  $\varpi^{4\cdot8^n}$. This proves that all powers of $X$ are
  distinct.
\end{proof}
Note that Georgi Genov and Plamen Siderov show
in~\cite{genov-siderov:grigorchuk} that $(1+A)(B+C)$, $(1+A)(B+D)$ and
$(1+A)(C+D)$ have infinite order in the group ring of $\Gg$. However, they
project to nil-elements in $\AA$.

Evidently $1+X$ belongs to the augmentation ideal $\varpi$, and is
also transcendental --- in particular, it is not nilpotent. However,
$\varpi$ contains many nilpotent elements:
\begin{prop}[\cite{vieira:modular}*{Theorem~2}]\label{prop:mononil}
  The semigroup $\{A,B,C,D\}^*\setminus\{1\}$ is nil of degree $8$.
\end{prop}
\begin{proof}
  Let $w\in S^n$ be a semigroup element. If $n$ is odd, then $w$ is
  either a $T\div T$ word or a $A\div A$ word, so $w^2=0$.

  If $n\equiv 2\pmod 4$, then either $w$ contains a $D$, in which case
  $w^2=0$ by Lemma~\ref{lem:alg:rel}, or $\psi(w^2)$ contains $D$'s in
  its non-zero entries, in which case $w^4=0$; or $\psi^2(w^4)$
  contains $D$'s in its non-zero entries, in which case $w^8=0$.

  Finally, if $n\equiv0\pmod 4$ and $n>0$, then $\psi(w)=\smallmat uuvv$
  or $\smallmat uvuv$, for some $u,v\in S^{n/2}$ with $uv=vu=0$. Then
  $\psi(w^8)=\smallmat{u^8}{u^8}{v^8}{v^8}$ or $\smallmat{u^8}{v^8}{u^8}{v^8}$,
  and we are done by induction on $n$.
\end{proof}

\subsubsection{Nillity and Primitivity of $\AA$}
To understand the representation theory of $\AA$, it is important to
determine whether $\AA$ is primitive. This depends on the Jacobson
radical of $\AA$, by the following simple result:
\begin{prop}\label{prop:sp=>p}
  If $\AA$ is semiprimitive, then it is primitive.
\end{prop}
\begin{proof}
  Since $\AA$ is semiprimitive, $\jac\AA=\bigcap_\PP\PP=0$, where the
  intersection is taken over all primitive ideals. However, if
  $\PP\neq0$ is primitive, then it has finite codimension
  by Theorem~\ref{thm:alg:ji}, so $\AA/\PP$ is finite-dimensional, and
  therefore nilpotent, because $\AA/\PP$ is the quotient of the group
  ring of a finite $2$-group, so $\PP=\varpi$. The only way to have
  $\jac\AA=0$ is therefore that $0$ be a primitive ideal.
\end{proof}

\begin{prop}\label{prop:=2:prim}
  If $\Bbbk$ is a field that is not algebraic over $\F[2]$, then
  $\AA$ is primitive.
\end{prop}
\begin{proof}
  Let $t$ be transcendental over $\F[2]$, and let
  $Y=A+B+AD\in\AA(\F[2])$ be transcendental, as in
  Theorem~\ref{thm:=2:laurent}. Assume for contradiction that
  $Y\in\jac\AA$. Then $1-tY$ is right invertible, i.e.\ there exists
  $r\in\AA$ with $(1-tY)r=1$. We may assume $r\in\AA(\F[2](t))$, so
  $(1-tY)p(t)=q(t)$ for $p(t)\in\AA(\F[2][t])$ and $q(t)\in\F[2][t]$.
  Again because $1-tY$ is invertible, we have
  $p(t)=q(t)\sum_{i=0}^\infty t^iY^i$. Considering this equality in
  degree higher than $\deg p$ and writing $q(t)=Q(t,1)$ as a
  homogeneous polynomial, we get $Q(1,Y)=0$ whence $Y$ is not
  transcendental.

  Therefore $\jac\AA\neq\varpi$, and $\jac\AA=0$ by
  Lemma~\ref{lem:jacalternative}, so $\AA$ is primitive by
  Proposition~\ref{prop:sp=>p}.
\end{proof}

Note that since $\AA$ is primitive for $\Bbbk=\F[2](t)$, it has a
maximal right ideal $L$ with trivial core, and therefore an
irreducible faithful nonprincipal\footnote{i.e.\ not of the form
  $e\AA$ for an idempotent $e$. Since $\AA$ is graded with
  $1$-dimensional degree-$0$ component, it has no idempotent except
  $0$ and $1$.} module $\AA/L$. One may take any maximal ideal $L$
containing $(1-tY)\AA$ with $Y$ as in the proof of
Proposition~\ref{prop:=2:prim}; however, there does not seem to be any
handy construction of such an $L$. On the other hand, the arguments
in~\cite{passman-t:reps}*{\S2} show that there are infinitely many
nonprincipal irreducible representations of $\AA$.

\begin{lem}[A. Smoktunowicz]\label{lem:jac<=>nil}
  Let $\JJ$ be a graded algebra (without unit) generated in degree
  $1$. Then the following are equivalent:
  \begin{enumerate}
  \item $\JJ$ is Jacobson radical;
  \item $M_n(\JJ)$ is graded nil\footnote{$M_n(\AA)$ is naturally
    $\Z^{n^2}$-graded by grading all entries independently.} for all
    $n$;
  \item $M_n(\JJ_1)$ is nil for all $n$, where $\JJ_1$ denotes the
  degree-$1$ component of $\JJ$.
  \end{enumerate}
\end{lem}
\begin{proof}
  We denote by $\FF$ the algebra $\JJ$ with a unit adjoined.  If $\JJ$
  is Jacobson radical, then $M_n(\JJ)$ is radical for all $n$. Take
  $x\in M_n(\FF)$, homogeneous of degree $d$. Then, since $1-x\in
  M_n(\FF)$ is invertible, the sum $\sum_{i\ge0}x^i$ must converge;
  now the component of degree $di$ of this sum is $x^i$; therefore
  $x^i=0$ for $i$ large enough.

  The next implication is obvious.

  Finally, assume $M_n(\JJ_1)$ is nil for all $n$, and choose
  $x\in\JJ$; write $x=x_1+\dots+x_r$ as a sum of
  monomials. Furthermore, write each monomial $x_i$, of degree $d_i$,
  as a product $x_i=x_{i,1}\dots x_{i,d_i}$ of monomials of degree
  $1$. Set
  \[A=\setsuch{(i,j)}{1\le i\le r,\,1\le j<d_i}\cup\{(0,0)\}.
  \]
  Construct the matrix $X\in M_A(\JJ_1)$ by
  \[X_{(i,j),(i',j')}=\begin{cases}
  x_{i,j+1} & \text{ if }i=i'\text{ and }j+1=j',\\
  x_{i',1} & \text{ if } (i,j)=(0,0)\text{ and }j'=1,\\
  x_{i,d_i} & \text{ if } j=d_i-1\text{ and }(i',j')=(0,0),\\
  \sum_{k:\,d_k=1}x_k & \text{ if }(i,j)=(i',j')=(0,0).
  \end{cases}
  \]
  Since $M_n(\JJ_1)$ is nil, there exists $N\in\N$ such that
  $X^N=0$. Now write formally $(1-x)^{-1}=1+y_1+y_2+\dots$ as a sum of
  homogeneous components. Then by induction
  $(X^s)_{(0,0),(i,j)}=y_{s-j}x_{i,1}\dots x_{i,j}$ if $i\ge1$ and
  $s>j$, and $y_s=(X^s)_{(0,0),(0,0)}$; therefore $y_s=0$ as soon as
  $s\ge N$, and $(1-x)^{-1}$ exists, so $x\in\jac\JJ$ and $\JJ$ is
  Jacobson radical.
\end{proof}

\begin{lem}\label{lem:Anil-MnAnil}
  The algebra $\AA$ is graded nil if and only if $M_n(\AA)$ is graded
  nil for all $n$.
\end{lem}
\begin{proof}
  Choose a homogeneous element $x\in M_n(\AA)$ of degree $\ge1$. It
  costs nothing to assume that $n$ is a power of two, say $n=2^t$.
  Then since $\varpi^3\le\KK$ by Proposition~\ref{prop:samefilt}, we
  have $x^3\in M_n(\KK)$, and therefore $y=\psi^{-t}(x^3)\in\KK$ is
  homogeneous. It follows that $y$, and therefore $x$, are nil
  elements.
\end{proof}

\begin{prop}\label{prop:nilprim}
  The algebra $\AA$ is non-primitive if and only if it is graded nil
  (i.e.\ all homogeneous elements of degree $\ge1$ are nil).
\end{prop}
\begin{proof}
  Assume first that $\AA$ is not graded nil. Then $\AA$ is not
  Jacobson radical by Lemma~\ref{lem:jac<=>nil}, so $\jac\AA=0$ by
  Lemma~\ref{lem:jacalternative}, and $\AA$ is primitive
  by~\ref{prop:sp=>p}.

  Assume next that $\AA$ is graded nil. Then $M_n(\AA)$ is graded nil
  by Lemma~\ref{lem:Anil-MnAnil}.  By Lemma~\ref{lem:jac<=>nil} the
  ideal $\varpi$ is Jacobson radical, so $\jac\AA=\varpi$.
\end{proof}

We denote below by $\AA_n$ the homogeneous part of $\AA$ of degree
$n$. The products $\AA_{n_1}\cdots\AA_{n_k}$ etc.\ are to be understood
as setwise products, and not linear spans of products.

\begin{lem}\label{lem:F2:nil}
  Let $n=n_1+\dots+n_k$ be even. Assume that for any choice of
  $n'_i,n''_i$ such that $|n_i-2n'_i|\le1$ and $|n_i-2n''_i|\le1$ and
  $\sum n'_i+n''_i=n$ we have
  \[(\AA_{n'_1}\AA_{n'_2}\cdots\AA_{n'_k}\AA_{n''_1}\AA_{n''_2}\cdots\AA_{n''_k})^t=0.\]
  Then
  \[(\AA_{n_1}\AA_{n_2}\cdots\AA_{n_k})^{2t+1}=0.\]
\end{lem}
\begin{proof}
  Choose $w_i\in\AA_{n_i}$, and write $w=w_1\cdots w_k$. For those
  $n_i$ which are even, we can write $w_i=w'_i+w''_i$ with $w'_i$ a
  linear combination of $A\div T$ words and $w''_i$ a linear
  combination of $T\div A$ words, while for the odd $n_i$ we can write
  $w_i=w'_i+w''_i$ with $w'_i$ a linear combination of $A\div A$ words
  and $w''_i$ a linear combination of $T\div T$ words.

  We next switch the $w'_i$ and $w''_i$ so that $w'_1$ is a $A\div$
  word, and $w'_{i+1}$ is a $A\div$ word if and only if $w'_i$ is a
  $\div T$ word. Then, since $n$ is even, $w^{2t}=0$ if and only if
  $(w'_1\cdots w'_n)^{2t}=0$ and $(w''_1\cdots w''_n)^{2t}=0$. We may
  therefore assume in turn that $w$ is a linear combination of $A\div
  T$ words, or is a linear combination of $T\div A$ words.

  We may also assume that each $w_i$ is either a linear combination of
  $A\div T$ words, or of $T\div A$ words, or of $A\div A$ words, or of
  $T\div T$ words.  We consider these cases in turn. If $w_i$ is a
  \begin{itemize}
  \item $A\div T$ word: then $\psi(w_i)=\smallmat uvuv=\vmat11\hmat
    uv$ and we set $x_i=u+v$;
  \item $T\div A$ word: then $\psi(w_i)=\smallmat uuvv=\vmat
    uv\hmat11$ and we set $x_i=u+v$;
  \item $A\div A$ word: then $\psi(w_i)=\smallmat
    uuuu=\vmat11u\hmat11$ and we set $x_i=u$;
  \item $T\div T$ word: then $\psi(w_i)=\smallmat uvwx$ and we set
    $x_i=u+v+w+x$.
  \end{itemize}

  If $w$ is a linear combination of $A\div T$ words, then
  $\psi(w^{2t})A=(x_1\cdots x_k)^{2t}A$, and by hypothesis
  $((x_1\cdots x_k)^2)^t=0$, so $w^{2t}A=0$.  If $w$ is a linear
  combination of $T\div A$ words, then $A\psi(w^{2t})=0$ by the same
  argument. In all cases $w^{2t+1}=0$.
\end{proof}

\begin{prop}\label{prop:F2:grnil}
  If $\Bbbk=\F[2]$, then $\AA$ is graded nil; more precisely, given
  $x\in\AA$ homogeneous of degree $n$, we have $x^{72n}=0$.
\end{prop}
\begin{proof}
  If $n$ is odd, we may replace $x$ by $x^2$, which will be of even
  degree $2n$. It is therefore sufficient to show that $x^{18n}=0$ for
  all homogeneous elements of even degree $n$, and from now on we
  assume that $n$ is even.

  Assume first that $x\in\AA_1^n$. Then $x=x_1\cdots x_n$, and since
  $\AA_1$ is spanned by $\{A,B,C,D\}$ with $B+C+D=0$, we may write
  $x_i=\alpha_i A+\beta_i T_i$ with $\alpha_i,\beta_i\in\F[2]$ and
  $T_i\in\{B,C,D\}$. Furthermore, since $n$ is even, we have $x=x'+x''$
  where $x'$ and $x''$ are monomials, with $x'\in(AT)^{n/2}$ and
  $x''\in(TA)^{n/2}$. Therefore $x^8=(x')^8+(x'')^8=0$ by
  Proposition~\ref{prop:mononil}.

  In the general case, set $t_0=8$ and $t_{i+1}=2t_i+1$. Then
  $t_i=9\cdot2^i-1$. Find $k\in\N$ such that $n\le2^k<2n$. Then,
  applying $k$ times Lemma~\ref{lem:F2:nil} to $x$, we have
  $x^{t_k}=0$, so \emph{a fortiori} $x^{18n}=0$.
\end{proof}

The following result answers a question in~\cite{vieira:modular}; it
also answers a conjecture attributed to
Goodearl~\cite{bell:examples}*{Conjecture~3.1}.
\begin{thm}\label{thm:primjacrad}
  If $\Bbbk$ is algebraic over $\F[2]$, then $\AA$ is graded nil and
  Jacobson radical.

  If $\Bbbk$ is not algebraic over $\F[2]$, then $\AA$ is not graded
  nil, and it is primitive.
\end{thm}
\begin{proof}
  Assume first that $\Bbbk$ is algebraic over $\F[2]$. Choose a
  homogeneous $x\in\AA(\Bbbk)$ of degree $\ge1$. Then
  $x\in\AA(\F[2^n])$ for some $n$, and therefore $x$ may be seen as a
  homogeneous element in $M_n(\AA)$, by embedding $\F[2^n]$ as a
  maximal field in $M_n(\F[2])$. Now $\AA$ is graded nil by
  Proposition~\ref{prop:F2:grnil}, so $M_n(\AA)$ is graded nil by
  Lemma~\ref{lem:Anil-MnAnil}, so $\AA(\F[2^n])$ is graded nil by
  restriction, and therefore $\varpi(\F[2^n])=\AA(\F[2^n])$ by
  Proposition~\ref{prop:nilprim} and Lemma~\ref{lem:jacalternative};
  finally $\AA(\Bbbk)$ is Jacobson radical since it is a union of such
  algebras.

  Assume now that $\Bbbk$ contains a transcendental element $t$. Then
  $\AA$ is primitive by Proposition~\ref{prop:=2:prim}, and the proof
  of Theorem~\ref{thm:=2:laurent}, just as
  Proposition~\ref{prop:nilprim}, imply that $\AA$ is not graded nil.
  Indeed the element $A+B+Dt$ has infinite order.
\end{proof}

\end{document}